\documentclass[11pt,leqno]{article}
\usepackage{amssymb}
\usepackage[mathscr]{eucal}
\usepackage{amsmath,amssymb,latexsym,theorem,bbm}
\usepackage{color}
\usepackage{appendix}

\newcommand{\CC}{\mathsf{C}}

\newcommand{\NN}{\mathbb{N}}

\newcommand{\RR}{\mathbb{R}}

\newcommand{\ZZ}{\mathbb{Z}}

\newcommand{\tP}{\widetilde{P}}

\newcommand{\bone}{{\boldsymbol{1}}}

\newcommand{\cA}{{\mathcal A}}
\newcommand{\cB}{{\mathcal B}}

\newcommand{\cC}{{\mathcal C}}

\newcommand{\cF}{{\mathcal F}}

\newcommand{\cL}{{\mathcal L}}

\newcommand{\dd}{\mathrm{d}}
\newcommand{\ee}{\mathrm{e}}

\newcommand{\EE}{\operatorname{\mathbb{E}}}
\newcommand{\PP}{\operatorname{\mathbb{P}}}

\newcommand{\tX}{\widetilde{X}}
\newcommand{\tY}{\widetilde{Y}}

\renewcommand{\mid}{\,|\,}

\renewcommand{\leq}{\leqslant}
\renewcommand{\geq}{\geqslant}

\newcommand{\distr}{\stackrel{\cL}{\longrightarrow}}

\newcommand{\proofend}{\hfill\mbox{$\Box$}}

\numberwithin{equation}{section}

\theoremstyle{change} \theorembodyfont{\em}
\newtheorem{Lem}{Lemma.}[section]
\newtheorem{Thm}[Lem]{Theorem.}

\theorembodyfont{\rm}
\newtheorem{Rem}[Lem]{Remark.}

\begin{document}

\begin{center}
 {\bfseries\Large
  Stationarity and ergodicity for an affine two factor model} \\[5mm]
 {\sc\large
  M\'aty\'as $\text{Barczy}^{*,\diamond}$,
  \ Leif $\text{D\"oring}^{**}$,
  \ Zenghu $\text{Li}^{***}$,
  \ Gyula $\text{Pap}^{****}$}
\end{center}

\vskip0.2cm

\noindent
 * Faculty of Informatics, University of Debrecen,
   Pf.~12, H--4010 Debrecen, Hungary.
   Tel.: +36-52-512900, Fax: +36-52-512996

\noindent
 ** Institut f\"ur Mathematik,
    Universit\"at Z\"urich,
    Winterthurerstrasse 190,
    CH-8057 Z\"urich,
    Switzerland.
    Tel.: +41-(0)44-63 55892, Fax: +41-(0)44-63 55705

\noindent
 *** School of Mathematical Sciences, Beijing Normal University,
     Beijing 100875,  People's Republic of China.
     Tel.: +86-10-58802900, Fax: +86-10-58808202

\noindent
 **** Bolyai Institute, University of Szeged,
      Aradi v\'ertan\'uk tere 1, H--6720 Szeged, Hungary.
      Tel.: +36-62-544033, Fax: +36-62-544548

\noindent e--mails: barczy.matyas@inf.unideb.hu (M. Barczy),
                    leif.doering@googlemail.com (L. D\"oring),\\
                    lizh@bnu.edu.cn (Z. Li),
                    papgy@math.u-szeged.hu (G. Pap).

\noindent $\diamond$ Corresponding author.

\renewcommand{\thefootnote}{}
\footnote{\textit{2010 Mathematics Subject Classifications\/}:
          60J25, 37A25.}
\footnote{\textit{Key words and phrases\/}:
 affine process, stationary distribution, ergodicity, Foster-Lyapunov criteria.}
\vspace*{0.2cm}
\footnote{The research of M. Barczy and G. Pap was realized in the frames of T\'AMOP 4.2.4.\ A/2-11-1-2012-0001 ,,National Excellence
Program -- Elaborating and operating an inland student and researcher personal support
system''. The project was subsidized by the European Union and co-financed by the
European Social Fund. 
Z. Li has been partially supported by NSFC under Grant No.\ 11131003 and 973
 Program under Grant No.\ 2011CB808001.}

\vspace*{-10mm}

\begin{abstract}
We study the existence of a unique stationary distribution and ergodicity for a \ $2$-dimensional affine process.
The first coordinate is supposed to be a so-called \ $\alpha$-root process with \ $\alpha\in(1,2]$.
The existence of a unique stationary distribution for the affine process is proved in case of \ $\alpha\in(1,2]$;
 \ further, in case of \ $\alpha=2$, \ the ergodicity is also shown.
\end{abstract}

\section{Introduction}

We consider the following 2-dimensional affine process (affine two factor model)
 \begin{align}\label{2dim_affine}
  \begin{cases}
   \dd Y_t = (a-bY_t)\,\dd t + \sqrt[\alpha]{Y_{t-}}\,\dd L_t,\qquad t\geq 0,\\
   \dd X_t = (m-\theta X_t)\,\dd t + \sqrt{Y_t}\,\dd B_t,\qquad t\geq 0,
  \end{cases}
 \end{align}
 where \ $a>0$, \ $b, \theta, m \in \RR$, \ $\alpha\in(1,2]$,
 \ $(L_t)_{t\geq 0}$ \ is a spectrally positive $\alpha$-stable L\'evy process
 with \ L\'evy measure \ $C_\alpha z^{-1-\alpha}\mathbf 1_{\{z>0\}}$ \
 with \ $C_\alpha:=(\alpha\Gamma(-\alpha))^{-1}$ \ (where \ $\Gamma$ \ denotes the Gamma function)
 in case \ $\alpha\in(1,2)$, \ a standard Wiener process in case \ $\alpha=2$, \ and
 \ $(B_t)_{t\geq 0}$ \ is an independent standard Wiener process.
Note that in case of \ $\alpha=2$, \ due to the almost sure continuity of the sample paths of a standard Wiener process,
 instead of \ $\sqrt{Y_{t-}}$ \ one can write \ $\sqrt{Y_t}$ \ in the first SDE of \eqref{2dim_affine}, and
 \ $Y$ \ is the so-called Cox-Ingersol-Ross (CIR) process; while in case of \ $\alpha\in(1,2)$, \
 \ $Y$ \ is called the \ $\alpha$-root process.
Note also that the process \ $(Y_t)_{t\geq 0}$ \ given by the first SDE of \eqref{2dim_affine} is a continuous state
 branching process with immigration with branching mechanism  \ $bz+\frac{1}{\alpha}z^\alpha$, $z\geq 0$, \
 and with immigration mechanism \ $az$, $z\geq 0$ \ (for more details, see the proof of Theorem \ref{Thm_ergodic1}, part (i)).
Chen and Joslin \cite{CheJos} have found several applications of the model \eqref{2dim_affine} with \ $\alpha=2$
 \ in financial mathematics, see their equations (25) and (26).

The process \ $(Y,X)$ \ given by \eqref{2dim_affine} is a special affine process.
The set of affine processes contains a large class of important Markov processes such as
 continuous state branching processes and Orstein-Uhlenbeck processes.
Further, a lot of models in financial mathematics are affine
 such as the Heston model \cite{Hes}, the model of Barndorff-Nielsen
 and Shephard \cite{BarShe} or the model due to Carr and Wu \cite{CarWu}.
A precise mathematical formulation and a complete characterization of regular
 affine processes are due to Duffie et al. \cite{DufFilSch}.
Later several authors have contributed to the theory of general affine processes:
 to name a few, Andersen and Piterbarg \cite{AndPit}, Dawson and Li \cite{DawLi},
 Filipovi\'{c} and Mayerhofer \cite{FilMay}, Glasserman and Kim \cite{GlaKim},
 Jena et al. \cite{JenKimXin} and Keller-Ressel et al. \cite{KelSchTei1}.

This article is devoted to study the existence of a unique stationary distribution and ergodicity of the affine process
 given by the SDE \eqref{2dim_affine}.
These kinds of results are important on their own rights, further they can be used for studying parameter estimation
 for the given model.
For the existing results on ergodicity of affine processes, see the beginning of Section \ref{Section_stat}.

Next we give a brief overview of the structure of the paper.
Section \ref{Section_preliminaires} is devoted to a preliminary discussion of
 the existence and uniqueness of a strong solution of the SDE \eqref{2dim_affine}
 by proving also that this solution is indeed an affine process, see, Theorem \ref{Pro_affine}.
In Section \ref{Section_stat} we prove the existence of a unique stationary distribution for the affine process
 given by \eqref{2dim_affine} in both cases \ $\alpha\in(1,2)$ \ and \ $\alpha=2$, \ provided that \ $a>0$,
 \ $b>0$ \ and \ $\theta>0$, \ see, Theorem \ref{Thm_ergodic1}.
In Section \ref{Section_ergod}, in case of \ $\alpha=2$, \ we prove ergodicity of the process in question
 provided that \ $a>0$, \ $b>0$ \ and \ $\theta>0$, \
 and we also show that the unique stationary distribution of the process is absolutely continuous, has finite (mixed)
 moments of any order by calculating some moments explicitly, too, see Theorems \ref{Thm_ergodic2} and \ref{Thm_ergodic3},
 respectively.

In a forthcoming paper we will use our results for studying parameter estimation for the process
 given by the SDE \eqref{2dim_affine}.

\section{The affine two factor model}\label{Section_preliminaires}

Let \ $\NN$, \ $\ZZ_+$, \ $\RR$ \ and \ $\RR_+$ \ denote the sets of positive integers, non-negative integers, real numbers and
 non-negative real numbers, respectively.
By \ $\|x\|$ \ and \ $\|A\|$ \ we denote the Euclidean norm of a vector \ $x\in\RR^m$ \ and the induced matrix norm
 \ $\Vert A\Vert=\sup\{\Vert Ax\Vert : x\in\RR^m, \; \Vert x\Vert=1\}$ \ of a matrix \ $A\in\RR^{n\times m}$, \ respectively.
By \ $\cC^2(\RR_+\times\RR, \RR)$, \ $\cC^2_c(\RR_+\times\RR, \RR)$ \ and \ $\cC^{\infty}_c(\RR_+\times\RR, \RR)$,
 \ we denote the set of twice continuously differentiable real-valued functions on \ $\RR_+\times\RR$, \
 the set of twice continuously differentiable real-valued functions on \ $\RR_+\times\RR$ \ with compact support
 and the set of infinitely differentiable real-valued functions on \ $\RR_+\times\RR$ \ with compact support, respectively.
Convergence in distribution will denoted by \ $\distr$.

Let \ $\big(\Omega, \cF, (\cF_t)_{t\geq 0}, \PP\big)$ \ be a filtered probability space
 satisfying the usual conditions, i.e.,  \ $(\Omega,\cF,\PP)$ \ is complete, the filtration \ $(\cF_t)_{t\geq 0}$ \ is right-continuous
 and \ $\cF_0$ \ contains all the $\PP$-null sets in \ $\cF$.
\ Let \ $(B_t)_{t\geq 0}$ \ be a standard \ $(\cF_t)_{t\geq 0}$-Wiener process and \ $(L_t)_{t\geq 0}$ \ be a spectrally positive
 \ $(\cF_t)_{t\geq 0}$-stable process with index \ $\alpha\in(1,2]$.
\ We assume that \ $B$ \ and \ $L$ \ are independent.
If \ $\alpha=2$, \ we understand that \ $L$ \ is a standard \ $(\cF_t)_{t\geq 0}$-Wiener process.
If \ $\alpha\in(1,2)$, \ we understand that \ $L$ \ is a \ $(\cF_t)_{t\geq 0}$-L\'evy process with L\'evy-Khintchine
 formula
 \[
  \EE(\ee^{iu L_1})
     = \exp\left\{ \int_0^{\infty}(\ee^{iuz} - 1 - iuz)C_\alpha z^{-1-\alpha}\,\dd z\right\},\quad u\in\RR,
 \]
 where \ $C_\alpha=(\alpha\Gamma(-\alpha))^{-1}$.
\ Recall that in case of \ $\alpha\in(1,2)$ \ the L\'evy-It\^{o} representation of \ $L$ \ takes the form
 \[
  L_t = \int_{(0,t]}\int_{(0,\infty)} z \widetilde{N}(\dd s,\dd z),\quad t\geq 0,
 \]
 where \ $\widetilde{N}(\dd s,\dd z)$ \ is a compensated Poisson random measure
 on \ $(0,\infty)^2$ \  with intensity measure
 \ $C_\alpha z^{-1-\alpha}\mathbf 1_{\{z>0\}}\,\dd s \,\dd z$.

\begin{Rem}\label{Rem_1}
We shed some light on the definition of the stochastic integral with respect
 to the spectrally positive $\alpha$-stable process \ $L$ \ in the first SDE of \eqref{2dim_affine}
 in case of \ $\alpha\in(1,2)$.
\ By Jacod and Shiryaev \cite[Corollary II.4.19]{JSh}, \ $L$ \ is a semimartingale so that
 Theorems I.4.31 and I.4.40 in Jacod and Shiryaev \cite{JSh} describe the classes of processes
 which are integrable with respect to \ $L$.
\ A more accessible integrability criteria is due to Kallenberg \cite[Theorem 3.1]{Kal1}.
Roughly speaking, a predictable process \ $V$ \ is locally integrable with respect to \ $L$ \
 (i.e., the stochastic integral \ $\int_0^t V_s\,\dd L_s$ \ exists for all \ $t\geq 0$) \
 if and only if \ $\int_0^t \vert V_s\vert^\alpha\,\dd s<\infty$ \ almost surely for all \ $t\geq 0$.
\ For the construction of stochastic integrals with respect to symmetric $\alpha$-stable processes, see
 also Rosinski and Woyczynski \cite[Theorem 2.1]{RosWoy}.
Another possible way is to consider the stochastic integral with respect to \ $L$ \ as a stochastic integral with respect to a certain
 compensated Poisson random measure, see the last equality on page 230 in Li \cite{Li}.
\proofend
\end{Rem}

The next proposition is about the existence and uniqueness of a strong solution of the SDE \eqref{2dim_affine}.

\begin{Thm}\label{Pro_affine}
Let \ $(\eta_0,\zeta_0)$ \ be a random vector independent of \ $(L_t,B_t)_{t\geq 0}$ \
 satisfying \ $\PP(\eta_0\geq 0) = 1$.
\ Then for all \ $a>0$, $b, m, \theta\in \RR$ \ and \ $\alpha\in(1,2]$, \ there is a (pathwise)
 unique strong solution \ $(Y_t,X_t)_{t\geq 0}$ \ of the SDE \eqref{2dim_affine} such that
 \ $\PP((Y_0,X_0) = (\eta_0,\zeta_0))=1$ \ and \ $\PP(Y_t\geq 0,\; \forall\;t\geq 0)=1$.
\ Further, we have
  \begin{align}\label{help50}
   Y_t = \ee^{-b(t-s)}\left( Y_s + a\int_s^t \ee^{-b(s-u)}\,\dd u
                                 + \int_s^t \ee^{-b(s-u)}\sqrt[\alpha]{Y_{u-}} \,\dd L_u\right)
 \end{align}
 for \ $0\leq s\leq t$, \ and
 \begin{align}\label{help2}
  X_t = \ee^{-\theta(t-s)}\left( X_s + m\int_s^t \ee^{-\theta(s-u)}\,\dd u
                                 + \int_s^t \ee^{-\theta(s-u)}\sqrt{Y_u} \,\dd B_u\right)
 \end{align}
 for \ $0\leq s\leq t$.
\ Moreover, \ $(Y_t,X_t)_{t\geq 0}$ \ is a regular affine process with infinitesimal generator
 \begin{align}\label{help76_infgen}
  \begin{split}
  (\cA f)(y,x)
   &= (a-by)f_1'(y,x) + (m-\theta x)f_2'(y,x)
      + \frac{1}{2}yf_{2,2}''(y,x) \\
   &\phantom{=\,}  + y \int_0^\infty \Big( f(y+z,x) - f(y,x) - z f_1'(y,x)\Big)
                              C_\alpha z^{-1-\alpha}\, \dd z
  \end{split}
 \end{align}
 in case of \ $\alpha\in(1,2)$, \ and
 \begin{align}\label{help77_infgen}
  (\cA f)(y,x)
   = (a-by)f_1'(y,x) + (m-\theta x)f_2'(y,x)
     + \frac{1}{2}y(f_{1,1}''(y,x) + f_{2,2}''(y,x))
 \end{align}
 in case of \ $\alpha=2$, \ where \ $(y,x)\in\RR_+\times\RR$, \ $f\in\cC^2_c(\RR_+\times\RR, \RR)$,
 \ and \ $f_i'$, \ $i=1,2$, \ and \ $f_{i,j}''$, \ $i,j\in\{1,2\}$,
 \ denote the first and second order partial derivatives of $f$ with respect to
 its \ $i$-th and \ $i$-th and \ $j$-th variables.
\end{Thm}

\noindent{\bf Proof.}
By Theorem 6.2 or Corollary 6.3 in Fu and Li \cite{FuLi} (for the case \ $\alpha\in(1,2)$) \ and
 by Yamada and Watanabe theorem (see, e.g., Karatzas and Shreve \cite[Proposition 5.2.13]{KarShr})
 and Ikeda and Watanabe \cite[Example 8.2, page 221]{IkeWat} (for the case \ $\alpha=2$),
 \ there is a (pathwise) unique non-negative strong solution \ $(Y_t)_{t\geq 0}$ \ of the first equation
 in \eqref{2dim_affine} with any initial value \ $\eta_0$ \ satisfying \ $\PP(\eta_0\geq 0)=1$.
\ An application of the It\^{o}'s formula to the process \ $(\ee^{b t} Y_t)_{t\geq 0}$ \ shows that
  \begin{align*}
  \dd(\ee^{b t} Y_t)
    &= b\ee^{b t}Y_t\,\dd t + \ee^{b t}\dd Y_t
     = b\ee^{b t}Y_t\,\dd t + \ee^{b t}\big( (a-b Y_t)\,\dd t + \sqrt[\alpha]{Y_{t-}}\,\dd L_t\big)\\
    &= a\ee^{bt}\,\dd t + \ee^{b t}\sqrt[\alpha]{Y_{t-}}\,\dd L_t,\quad t\geq 0.
 \end{align*}
Hence, if \ $0\leq s\leq t$, \ then
  \[
    \ee^{b t} Y_t
      = \ee^{b s} Y_s + a\int_s^t \ee^{b u}\,\dd u
                                 + \int_s^t \ee^{b u}\sqrt[\alpha]{Y_{u-}} \,\dd L_u,
  \]
 and \eqref{help50} follows.

Let \ $(X_t)_{t\ge 0}$ \ be defined by \eqref{help2}.
By It\^{o}'s formula, we have
 \begin{align*}
 \dd X_t
      & = -\theta \ee^{-\theta t}\Big(X_0 + m\int_0^t \ee^{\theta s}\,\dd s
                                    + \int_0^t \ee^{\theta s}\sqrt{Y_s}\,\dd B_s\Big)\,\dd t
         + \ee^{-\theta t} (m\ee^{\theta t} \,\dd t + \ee^{\theta t}\sqrt{Y_t}\,\dd B_t)  \\
      & = -\theta \ee^{-\theta t}\Big(X_0 + m\int_0^t \ee^{\theta s}\,\dd s
           + \int_0^t \ee^{\theta s}\sqrt{Y_s}\,\dd B_s\Big)\dd t + m\,\dd t + \sqrt{Y_t}\,\dd B_t, \quad t\geq 0.
 \end{align*}
Then the process \ $(X_t)_{t\ge 0}$ \ is indeed a strong solution of the second equation in \eqref{2dim_affine},
 giving the existence of the strong solution of \eqref{2dim_affine}.
Conversely, if \ $(X_t)_{t\geq 0}$ \ is a strong solution of the second equation in \eqref{2dim_affine},
 by an application of It\^o's formula to the process \ $(\ee^{\theta t}X_t)_{t\geq 0}$, \ we have
 \begin{align*}
  \dd(\ee^{\theta t} X_t)
    &= \theta\ee^{\theta t}X_t\,\dd t + \ee^{\theta t}\dd X_t
     = \theta\ee^{\theta t}X_t\,\dd t + \ee^{\theta t}\big( (m-\theta X_t)\,\dd t + \sqrt{Y_t}\,\dd B_t\big)\\
    &= m\ee^{\theta t}\,\dd t + \ee^{\theta t}\sqrt{Y_t}\,\dd B_t,\quad t\geq 0.
 \end{align*}
If \ $0\leq s\leq t$, \ then
 \[
   \ee^{\theta t} X_t
      = \ee^{\theta s} X_s + m\int_s^t \ee^{\theta u}\,\dd u
                                 + \int_s^t \ee^{\theta u}\sqrt{Y_u} \,\dd B_u,
 \]
 and hence \ $(X_t)_{t\ge 0}$ \ must be given by \eqref{help2}.
This proves the pathwise uniqueness for the second equation in \eqref{2dim_affine}, and hence
 the pathwise uniqueness for \eqref{2dim_affine}.

Now we turn to check that \ $(Y_t,X_t)_{t\geq 0}$ \ is a regular affine process with the given infinitesimal generator.
We may and do suppose that the initial value is deterministic, say,
 \ $(Y_0, X_0) = (y_0, x_0) \in \RR_+ \times \RR$, \ since the infinitesimal
 generator of a time homogeneous Markov process does not depend on the initial
 value of the Markov process.

In case of \ $\alpha=2$, \ by It\^{o}'s formula, for all \ $f\in\cC^2_c(\RR_+\times\RR, \RR)$ \ we have
 \begin{align*}
  f(Y_t,X_t) & = f(y_0, x_0)
                + \int_0^t f_1'(Y_s,X_s)\sqrt{Y_s} \,\dd L_s
                + \int_0^t f_2'(Y_s,X_s)\sqrt{Y_s} \,\dd B_s \\
            &\phantom{=\;} + \int_0^t f_1'(Y_s,X_s)(a-bY_s)\,\dd s
                 + \int_0^t f_2'(Y_s,X_s)(m-\theta X_s)\,\dd s \\
            &\phantom{=\;} + \frac{1}{2}\left(\int_0^t f_{1,1}''(Y_s,X_s)Y_s\,\dd s
                                   + \int_0^t f_{2,2}''(Y_s,X_s)Y_s\,\dd s \right) 
 \end{align*}
 \begin{align*}                                  
           = f(y_0, x_0) + \int_0^t (\cA f)(Y_s,X_s)\,\dd s
             + M_t(f),\qquad t\geq 0,
 \end{align*}
 where
 \[
    M_t(f):= \int_0^t f_1'(Y_s,X_s)\sqrt{Y_s} \,\dd L_s
               + \int_0^t f_2'(Y_s,X_s)\sqrt{Y_s} \,\dd B_s,\qquad t\geq 0,
 \]
 and \ $\cA f$ \ is given by \eqref{help77_infgen}.
It remains to show that \ $(M_t(f))_{t\geq 0}$ \ is a local martingale with
 respect to the filtration \ $(\cF_t)_{t\geq 0}$.
\ However, we prove that it is a square integrable martingale
 with respect to the filtration \ $(\cF_t)_{t\geq 0}$, \ since
 \begin{align*}
   &\int_0^t \EE((f_1'(Y_s,X_s))^2Y_s) \,\dd s
       \leq C_1    \int_0^t \EE(Y_s) \,\dd s
       <\infty,\qquad t\geq 0,\\
   &   \int_0^t \EE((f_2'(Y_s,X_s))^2Y_s) \,\dd s
       \leq C_2    \int_0^t \EE(Y_s) \,\dd s
       <\infty,\qquad t\geq 0,
 \end{align*}
 with some constants \ $C_1>0$ \ and \ $C_2>0$, \ where the finiteness of the integrals follows by
 that
 \[
   \EE(Y_s) = \ee^{-bs} y_0 + a\int_0^s \ee^{-bu}\,\dd u, \qquad s\geq 0,
 \]
 see, e.g., Cox et al. \cite[Equation (19)]{CoxIngRos}, Jeanblanc et al. \cite[Theorem 6.3.3.1]{JeaYorChe}
 or Proposition 3.2 in Barczy et al.~\cite{BarDorLiPap}.

In case of \ $\alpha\in(1,2)$, \ by the L\'evy-It\^{o} representation of \ $L$, \ we can rewrite the SDE \eqref{2dim_affine}
 into the integral form
 \begin{align*}
  \begin{cases}
   Y_t = y_0 + \int_0^t (a-bY_s)\,\dd s + \int_0^t\int_0^\infty z\sqrt[\alpha]{Y_{s-}}\widetilde N (\dd s,\dd z),\qquad t\geq 0,\\
   X_t = x_0 + \int_0^t (m-\theta X_s)\,\dd s + \int_0^t \sqrt{Y_s}\,\dd B_s,\qquad t\geq 0.
  \end{cases}
 \end{align*}
By It\^{o}'s formula, for all \ $f\in\cC^2_c(\RR_+\times\RR, \RR)$ \ we have
 \begin{align*}
  &f(Y_t,X_t)
     = f(y_0, x_0)
       + \int_0^t f_1'(Y_s,X_s)(a-bY_s)\,\dd s
       + \int_0^t f_2'(Y_s,X_s)(m-\theta X_s)\,\dd s \\
    &\phantom{=\;}
       + \frac{1}{2}\int_0^t f_{2,2}''(Y_s,X_s) Y_s\,\dd s
       + \int_0^t f_2'(Y_s,X_s)\sqrt{Y_s} \,\dd B_s \\
   &\phantom{=\;}
       + \int_0^t\int_0^\infty \Big( f(Y_{s-}+z\sqrt[\alpha]{Y_{s-}},X_{s-}) - f(Y_{s-},X_{s-}) \Big)
                              \,\widetilde N(\dd s,\dd z)\\
    &\phantom{=\;}
      + \int_0^t\int_0^\infty \Big( f(Y_s+z\sqrt[\alpha]{Y_s},X_s) - f(Y_s,X_s)
                                     - z\sqrt[\alpha]{Y_s}f_1'(Y_s,X_s) \Big)
                                     C_\alpha z^{-1-\alpha}\,\dd s\dd z\\
    &= f(y_0, x_0) + \int_0^t (\cA f)(Y_s,X_s)\,\dd s + M_t(f),\qquad t\geq 0,
 \end{align*}
 where
 \begin{align*}
   M_t(f)&:= \int_0^t f_2'(Y_s,X_s)\sqrt{Y_s} \,\dd B_s\\
         &\phantom{:=\;} + \int_0^t\int_0^\infty \Big( f(Y_{s-}+z\sqrt[\alpha]{Y_{s-}},X_{s-}) - f(Y_{s-},X_{s-}) \Big)
                              \,\widetilde N(\dd s,\dd z),\qquad t\geq 0,
 \end{align*}
 and, by the change of variable
 \ $z\sqrt[\alpha]{y}:=\widetilde z$, \
 \begin{align*}
  (\cA f)(y,x)
   &:= (a-by)f_1'(y,x) + (m-\theta x)f_2'(y,x)
      + \frac{1}{2}yf_{2,2}''(y,x) \\
   &\phantom{:=\;}  + y\int_0^\infty \Big( f(y+\widetilde z,x) - f(y,x) - \widetilde z f_1'(y,x)\Big)
                              C_\alpha z^{-1-\alpha}\,\dd \widetilde z
 \end{align*}
 for \ $(y,x)\in\RR_+\times\RR$ \ and \ $f\in\cC^2_c(\RR_+\times\RR, \RR)$.
\ It remains to show that \ $(M_t(f))_{t\geq 0}$ \ is a local martingale
 with respect to the filtration \ $(\cF_t)_{t\geq 0}$.
\ However, we prove that it is
 a martingale with respect to the filtration \ $(\cF_t)_{t\geq 0}$.
\ The martingale property of \ $\int_0^t f_2'(Y_s,X_s)\sqrt{Y_s} \,\dd B_s$, $t\geq 0$, \ follows in the very
 same way as in the case of \ $\alpha=2$ \ using that there exists some constant \ $C_3>0$ \ such that
 \ $\EE(Y_t)\leq C_3(1+ y_0\ee^{-bt/\alpha})$ \ for all \ $t\geq 0$, \ see, e.g., formula (2.7) or
 Proposition 2.10 in Li and Ma \cite{LiMa}.
Now we turn to check that
 \[
  M^1_t(f):=\int_0^t\int_0^\infty \Big( f(Y_{s-}+z\sqrt[\alpha]{Y_{s-}},X_{s-}) - f(Y_{s-},X_{s-}) \Big)
                              \,\widetilde N(\dd s,\dd z),\quad t\geq 0,
 \]
 is a martingale.
For all \ $n\in\NN$, \ let
 \begin{align*}
  &M^{2,n}_t(f) := \int_0^t\int_1^\infty \Big( f\big((Y_{s-}\wedge n)+z\sqrt[\alpha]{Y_{s-}\wedge n},X_{s-}\big) - f(Y_{s-}\wedge n,X_{s-}) \Big)
                              \,\widetilde N(\dd s,\dd z),\\
  &M^{3,n}_t(f) := \int_0^t\int_0^1 \Big( f\big((Y_{s-}\wedge n)+z\sqrt[\alpha]{Y_{s-}\wedge n},X_{s-}\big) - f(Y_{s-}\wedge n,X_{s-}) \Big)
                              \,\widetilde N(\dd s,\dd z).
 \end{align*}
By Taylor's theorem, we have
 \begin{align*}
  |f\big((Y_{s-}\wedge n)+z\sqrt[\alpha]{Y_{s-}\wedge n},X_{s-}\big) - f(Y_{s-}\wedge n,X_{s-})|
  \leq z \sqrt[\alpha]{Y_{s-}\wedge n} \sup_{(y,x)\in\RR_+\times\RR} |f_1'(y,x)|
 \end{align*}
 for \ $z \in \RR_+$.
\ Since
 \begin{align*}
   &\EE\left( \int_0^t\int_1^\infty
               \Big| f\big((Y_{s-}\wedge n)+z\sqrt[\alpha]{Y_{s-}\wedge n},X_{s-}\big) - f(Y_{s-}\wedge n,X_{s-}) \Big|
                              \,C_\alpha z^{-1-\alpha}\dd s\,\dd z\right)\\
   &\quad\leq \biggl(\sup_{(y,x)\in\RR_+\times\RR} |f_1'(y,x)|\biggr)
              \int_0^t\int_1^\infty
               \EE\left( \sqrt[\alpha]{Y_s\wedge n}\right)
                C_\alpha z^{-\alpha}\dd s\,\dd z\\
   &\quad\leq C_4 n^{1/\alpha}t \int_1^\infty z^{-\alpha}\,\dd z
     = C_4 n^{1/\alpha}\frac{t}{\alpha-1} <\infty,\quad t\in\RR_+,
 \end{align*}
 with some constant \ $C_4>0$, \ by Lemma 3.1 in Chapter II and page 62 in Ikeda and Watanabe \cite{IkeWat},
 we get \ $(M^{2,n}_t(f))_{t\geq 0}$ \ is a martingale with respect to the filtration \ $(\cF_t)_{t\geq 0}$.
\ Further, since
 \begin{align*}
  &\EE\left( \int_0^t\int_0^1 \Big( f\big((Y_{s-}\wedge n)+z\sqrt[\alpha]{Y_{s-}\wedge n},X_{s-}\big) - f(Y_{s-}\wedge n,X_{s-}) \Big)^2
                              \,C_\alpha z^{-1-\alpha}\dd s\,\dd z \right)\\
  &\quad \leq C_5 \int_0^t \EE((Y_s\wedge n)^{2/\alpha})\,\dd s
      \int_0^1 z^{1-\alpha}\,\dd z
   \leq C_5 n^{2/\alpha} \frac{t}{2-\alpha}<\infty,\quad t\in\RR_+,
 \end{align*}
 with some constant \ $C_5>0$, \ by pages 62 and 63 in Ikeda and Watanabe \cite{IkeWat},
 we get \ $(M^{3,n}_t(f))_{t\geq 0}$ \ is a square integrable martingale with respect to the filtration \ $(\cF_t)_{t\geq 0}$.
This yields the martingale property of \ $(M^1_t(f))_{t\geq 0}$ \ with respect to the filtration
 \ $(\cF_t)_{t\geq 0}$.
\ Indeed, for all \ $n\in\NN$ \ and \ $f\in\cC_c^2(\RR_+\times\RR,\RR)$, \ let
 \begin{align*}
  \eta_t^n(f) := M^1_t(f) - M^{2,n}_t(f) - M^{3,n}_t(f),\quad t\in\RR_+.
 \end{align*}
Then, by Taylor's theorem,
 \begin{align*}
   \eta_t^n(f)
    &= \int_0^t\int_0^\infty \bone_{\{Y_{s-}>n\}}\Big( f(Y_{s-}+z\sqrt[\alpha]{Y_{s-}},X_{s-}) - f(Y_{s-},X_{s-}) \Big)
                              \,\widetilde N(\dd s,\dd z)\\
    &= \int_0^t\int_0^\infty \bone_{\{Y_{s-}>n\}}
                                   f_1'\big(Y_{s-}+\zeta z\sqrt[\alpha]{Y_{s-}},X_{s-}\big)
                                   z \sqrt[\alpha]{Y_{s-}}\,\widetilde N(\dd s,\dd z) \\
    &= \int_0^t
                 f_1'\big(Y_{s-}+\zeta z\sqrt[\alpha]{Y_{s-}},X_{s-}\big)
                 \sqrt[\alpha]{Y_{s-}} \,\dd L_s,\qquad t\in\RR_+
 \end{align*}
 with some (not necessarily measurable) function \ $\zeta:\Omega\to\RR$.
\ Note that, despite of the fact that \ $\zeta$ \ is not necessarily
 measurable, \ $f_1'\big(Y_{s-}+\zeta z\sqrt[\alpha]{Y_{s-}},X_{s-}\big)$ \ is a
 random variable, since it equals
 \ $f(Y_{s-}+z\sqrt[\alpha]{Y_{s-}},X_{s-}) - f(Y_{s-},X_{s-})$,
 \ which is a random variable.
In the same way,
 \ $\big(f_1'(Y_{s-}+\zeta z\sqrt[\alpha]{Y_{s-}},X_{s-})\big)_{s\in[0,t]}$
 \ is a predictable process, thus one can use Lemma 2.8 in Li and
 Ma \cite{LiMa}, hence there exist some constants \ $C_6>0$ \ and \ $C_7>0$ \ such that
 for all \ $t\in\RR_+$,
 \begin{align*}
  \EE\left( \sup_{s\in[0,t]}\vert \eta_s^n(f)\vert \right)
  &\leq C_6 \EE\left( \left( \int_0^t\bone_{\{Y_s>n\}} \big(f_1'\big(Y_s+\zeta z\sqrt[\alpha]{Y_s},X_s\big)\big)^\alpha Y_s \,\dd s\right)^{1/\alpha}\right)\\
  &\leq C_7 \left(\int_0^t \EE\left(\bone_{\{Y_s>n\}} Y_s\right)\dd s \right)^{1/\alpha},
 \end{align*}
 where the last inequality follows by Jensen inequality using also that \ $f'$ \ is bounded.
Using that there exists some constant \ $C_8>0$ \ such that
 \ $\EE(Y_t) \leq C_8 (1 + y_0 \ee^{-bt/\alpha})$, \ $t\in\RR_+$ \ (see, e.g.,
 formula (2.7) or Proposition 2.10 in Li and Ma \cite{LiMa}), we have
 \[
    \int_0^t \EE\left(\bone_{\{Y_s>n\}} Y_s\right) \,\dd s
       \leq C_8 \int_0^t (1 + y_0 \ee^{-bs/\alpha})\,\dd s<\infty,\qquad n\in\NN,
 \]
 and then, by dominated convergence theorem,
 \begin{align}\label{help78}
  \lim_{n\to\infty}  \EE\left( \sup_{s\in[0,t]}\vert \eta_s^n(f)\vert \right) = 0,
  \qquad t\in\RR_+.
 \end{align}
This yields that \ $(M^1_t(f))_{t\geq 0}$ \ is a martingale with respect to the filtration \ $(\cF_t)_{t\geq 0}$.
 \ Indeed, for all \ $0\leq s\leq t$ \ and \ $A\in\cF_s$, \ by \eqref{help78}, we have
 \begin{align*}
  &\EE(\vert \eta^n_t(f)\bone_A\vert) \leq \EE(\vert \eta^n_t(f)\vert) \to 0 \quad \text{as \ $n\to\infty$,}\\
  &\EE(\vert \eta^n_s(f)\bone_A\vert) \leq \EE(\vert \eta^n_s(f)\vert) \to 0 \quad \text{as \ $n\to\infty$.}
 \end{align*}
Hence
 \begin{align}\label{help79}
  \begin{split}
  &\lim_{n\to\infty} \EE((M^{2,n}_t(f) + M^{3,n}_t(f))\bone_A)
                      =  \EE(M^1_t(f)\bone_A), \\
  &\lim_{n\to\infty} \EE((M^{2,n}_s(f) + M^{3,n}_s(f))\bone_A)
                      =  \EE(M^1_s(f)\bone_A).
 \end{split}
 \end{align}
Since \ $(M^{2,n}_t(f) + M^{3,n}_t(f))_{t\geq 0}$ \ is a martingale with respect to the filtration
 \ $(\cF_t)_{t\geq 0}$ \ for all \ $n\in\NN$, \ we have
 \[
  \EE\big( (M^{2,n}_t(f) + M^{3,n}_t(f))\bone_A\big)
     = \EE\big( (M^{2,n}_s(f) + M^{3,n}_s(f))\bone_A\big), \quad n\in\NN,
 \]
 which, by \eqref{help79}, yields that
 \[
   \EE\big(M^1_t(f)\bone_A\big)
     = \EE\big( M^1_s(f)\bone_A\big)
 \]
 for all \ $0\leq s\leq t$ \ and \ $A\in\cF_s$, \ i.e., \ $(M^1_t(f))_{t\geq 0}$ \ is a martingale with respect to the filtration
 \ $(\cF_t)_{t\geq 0}$.

Finally, we check that the transition semigroup \ $(P_t)_{t\geq 0}$ \ with state space \ $\RR_+\times\RR$ \ corresponding
 to \ $(Y_t,X_t)_{t\geq 0}$ \ is a regular affine semigroup having infinitesimal generator given by
 \eqref{help76_infgen} and  \eqref{help77_infgen} according to the cases \ $\alpha\in(1,2)$ \ and
 \ $\alpha=2$.
With the notations of Dawson and Li \cite{DawLi},
 \ $(0,(\alpha_{i,j})_{i,j=1}^2,(b_i)_{i=1}^2,(\beta_{i,j})_{i,j=1}^2,0,\mu)$ \ and
 \ $(0,(\alpha_{i,j})_{i,j=1}^2,(b_i)_{i=1}^2,(\beta_{i,j})_{i,j=1}^2,0,0)$ \
 are sets of admissible parameters according to the cases \ $\alpha\in(1,2)$ \ and \ $\alpha=2$,
 \ where
 \begin{align*}
   &\begin{pmatrix}
    \alpha_{1,1} & \alpha_{1,2} \\
    \alpha_{2,1} & \alpha_{2,2} \\
  \end{pmatrix}
  :=\begin{cases}
     \begin{pmatrix}
      0 & 0 \\
      0 & \frac{1}{2} \\
     \end{pmatrix}
       & \text{if \ $\alpha\in(1,2)$,}\\[2mm]
    \begin{pmatrix}
      \frac{1}{2} & 0 \\
       0 & \frac{1}{2} \\
     \end{pmatrix}
       & \text{if \ $\alpha=2$,}
   \end{cases}\\[1mm]
 & (b_1,b_2) := (a,m) \qquad \text{for \ $\alpha\in(1,2]$,}\\[1mm]
 & \begin{pmatrix}
    \beta_{1,1} & \beta_{1,2} \\
    \beta_{2,1} & \beta_{2,2} \\
  \end{pmatrix}
  :=
   \begin{pmatrix}
    b & 0 \\
    0 & -\theta \\
  \end{pmatrix}
  \qquad \text{for \ $\alpha\in(1,2]$,}\\[1mm]
 & \mu(\dd y,\dd x) := C_\alpha y^{-1-\alpha}\dd y\delta_0(\dd x)
   \qquad \text{ for \ $\alpha\in(1,2)$,}
  \end{align*}
 where \ $\delta_0$ \ denotes the Dirac measure concentrated on \ $0\in\RR$.
\ Indeed, condition (vi) of Definition 6.1 in Dawson and Li \cite{DawLi} holds, since
 \begin{align*}
  \int_0^\infty \int_{-\infty}^\infty
    & (\vert y\vert\wedge y^2) \mu(\dd y,\dd x)
     + \int_0^\infty \int_{-\infty}^\infty
     (\vert x\vert\wedge x^2) \mu(\dd y,\dd x) \\
  & = C_\alpha \int_0^\infty \int_{-\infty}^\infty (\vert y\vert\wedge y^2) y^{-1-\alpha}
     \,\dd y \delta_0(\dd x)
     = C_\alpha \int_0^\infty (y\wedge y^2)\,y^{-1-\alpha}\,\dd y \\
  & = C_\alpha \int_0^1 y^{1-\alpha}\,\dd y
     + C_\alpha \int_1^\infty y^{-\alpha}\,\dd y
    <\infty.
 \end{align*}
Hence Theorem 2.7 in Duffie et al.~\cite{DufFilSch} (see also Theorem 6.1 in Dawson and Li \cite{DawLi})
 yields that for these sets of admissible parameters, there exists a regular affine semigroup
 \ $(Q_t)_{t\geq 0}$ \ with infinitesimal generator given by \eqref{help76_infgen} and \eqref{help77_infgen}
 according to the cases \ $\alpha\in(1,2)$ \ and \ $\alpha=2$.
\ By Theorem 2.7 in Duffie et al.~\cite{DufFilSch}, \ $\cC^\infty_c(\RR_+\times\RR, \RR)$ \ is a core
 of the infinitesimal generator corresponding to the affine semigroup \ $(Q_t)_{t\geq 0}$.
Since we have checked that the infinitesimal generators corresponding to
 the transition semigroups \ $(P_t)_{t\geq 0}$ \ and \ $(Q_t)_{t\geq 0}$ \
 (defined on the Banach space of bounded functions on \ $\RR_+\times\RR$)
 \ coincide on \ $\cC^\infty_c(\RR_+\times\RR, \RR)$, \ by the definition of a core, we get they coincide
 on the Banach space of bounded functions on \ $\RR_+\times\RR$.
\ This yields that \ $(Y_t,X_t)_{t\geq 0}$ \ is an affine process with infinitesimal generator \eqref{help76_infgen}
 and \eqref{help77_infgen} according to the cases \ $\alpha\in(1,2)$ \ and \ $\alpha=2$.
We also note that we could have used Lemma 10.2 in Duffie et al.~\cite{DufFilSch} for concluding that
 \ $(Y_t,X_t)_{t\geq 0}$ \ is a regular affine process with infinitesimal generator \eqref{help76_infgen}
 and \eqref{help77_infgen} according to the cases \ $\alpha\in(1,2)$ \ and \ $\alpha=2$, \ since we have checked
 that \ $(M_t(f))_{t\geq 0}$ \ is a martingale with respect
 to the filtration \ $(\cF_t)_{t\geq 0}$ \ for any \ $f\in\cC^2_c(\RR_+\times\RR, \RR)$ \
 in both cases \ $\alpha\in(1,2)$ \ and \ $\alpha=2$.
\proofend

\begin{Rem}
Note that in Theorem \ref{Pro_affine} it is the assumption \ $a>0$ \ which
 ensures \ $\PP(Y_t\geq 0,\; \forall\;t\geq 0)=1$.
\proofend
\end{Rem}

\section{Stationarity}\label{Section_stat}

The study of existence of stationary distributions for affine processes in general is currently under active research.

In the special case of continuous state branching processes with immigration the question of existence of
 a unique stationary distribution has been well-studied, see Li \cite[Theorem 3.20 and  Corollary 3.21]{Li}
 or Keller-Ressel and Mijatovi\'{c} \cite[Theorem 2.6]{KelMij}.

Glasserman and Kim \cite[Theorem 2.4]{GlaKim} proved existence of a unique stationary distribution for the process
  \begin{align}\label{2dim_affine_GlaKim}
  \begin{cases}
   \dd Y_t = (a-bY_t)\,\dd t + \sqrt{Y_t}\,\dd L_t,\qquad t\geq 0,\\
   \dd X_t = -\theta X_t\,\dd t + \sqrt{1+\sigma Y_t}\,\dd B_t,\quad t\geq 0,
  \end{cases}
 \end{align}
 where $a>0$, $b>0$, $\theta>0$, $\sigma\geq 0$ \ and \ $L$ \ and \ $B$ \ are independent standard Wiener processes.

The following result states the existence of a unique stationary distribution of the affine process given by
 the SDE \eqref{2dim_affine} for both cases \ $\alpha\in(1,2)$ \ and \ $\alpha=2$.

\begin{Thm}\label{Thm_ergodic1}
Let us consider the 2-dimensional affine model \eqref{2dim_affine} with
 \ $a>0$, \ $b>0$, $m\in\RR$, \ $\theta>0$, \ and with a random initial value
 \ $(\eta_0,\zeta_0)$ \ independent of \ $(L_t,B_t)_{t\geq 0}$
 \ satisfying \ $\PP(\eta_0\geq 0) = 1$.
\ Then
 \renewcommand{\labelenumi}{{\rm(\roman{enumi})}}
 \begin{enumerate}
 \item \ $(Y_t,X_t)\distr (Y_\infty,X_\infty)$ \ as \ $t\to\infty$, \ and the distribution of
       \ $(Y_\infty,X_\infty)$ \ is given by
       \begin{align}\label{help18}
          \EE\big(\ee^{-\lambda_1 Y_\infty + i\lambda_2 X_\infty}\big)
            = \exp\left\{-a\int_0^\infty v_s(\lambda_1,\lambda_2)\,\dd s + i\frac{m}{\theta}\lambda_2\right\}
       \end{align}
       for \ $(\lambda_1,\lambda_2)\in\RR_+\times\RR$, \ where \ $v_t(\lambda_1,\lambda_2)$, $t\geq 0$, \ is the unique
       non-negative solution of the (deterministic) differential
       equation
       \begin{align}\label{DE1}
         \begin{cases}
          \frac{\partial v_t}{\partial t} (\lambda_1,\lambda_2)
               = -bv_t(\lambda_1,\lambda_2) - \frac{1}{\alpha}(v_t(\lambda_1,\lambda_2))^\alpha
                  + \frac{1}{2}\ee^{-2\theta t}\lambda_2^2,\qquad t\geq 0,\\
          v_0(\lambda_1,\lambda_2) = \lambda_1.
         \end{cases}
       \end{align}
 \item supposing that the random initial value \ $(\eta_0,\zeta_0)$ \ has the same distribution as \ $(Y_\infty,X_\infty)$
       \ given in part (i), we have \ $(Y_t,X_t)_{t\geq 0}$ \ is strictly stationary.
\end{enumerate}
\end{Thm}

\noindent{\bf Proof.}
{\bf (i):} {\bf Step 1.}
In this step we give some representations of the affine transition semigroup
 \ $(P_t)_{t\geq 0}$ \ with state space \ $\RR_+\times\RR$ \ corresponding to the
 process given by the SDE \eqref{2dim_affine}.
By Theorem 6.1 in Dawson and Li \cite{DawLi} and Theorem \ref{Pro_affine}, the transition
 semigroup \ $(P_t)_{t\geq 0}$ \ is given by
 \begin{align}\label{help12}
  \int_{\RR_+\times\RR} \ee^{\langle u,\xi\rangle}P_t((y_0,x_0),\dd\xi)
     = \ee^{\langle (y_0,x_0),\psi_t(u)\rangle + \phi_t(u)}
 \end{align}
 for \ $u\in U$, $(y_0,x_0)\in\RR_+\times\RR,\;t\geq 0$, \ where \ $U:=\CC_{-}\times(i\RR)$ \ with
 \begin{align*}
  \CC_{-} := \{z_1+iz_2 : z_1\in(-\infty,0],\; z_2\in\RR\},\qquad i\RR :=\{iz_2 : z_2\in\RR\},
 \end{align*}
 and for all \ $u=(u_1,u_2)\in U$, \ we have \ $\psi_t(u)=(\psi^{(1)}_t(u),\ee^{-\theta t}u_2)$, $t\geq 0$, \
 where \ $\psi^{(1)}_t(u)$, $t\geq 0$, \ is a solution of the generalized Riccati
 equation
 \begin{align}\label{help11}
  \begin{cases}
   \frac{\partial \psi^{(1)}_t}{\partial t}(u) = R(\psi^{(1)}_t(u),\ee^{-\theta t}u_2),\qquad t\geq 0,\\
   \psi^{(1)}_0(u) = u_1,
  \end{cases}
 \end{align}
 and
 \begin{align*}
  \phi_t(u) = \int_0^t F(\psi_s^{(1)}(u),\ee^{-\theta s}u_2)\,\dd s, \qquad t\geq 0.
 \end{align*}
Here for \ $\alpha\in(1,2]$, \ the (complex valued) functions \ $F$ \ and \ $R$ \ are given by
  \begin{align*}
  F(u) = au_1 + mu_2, \qquad R(u) = -bu_1 + \frac{(-u_1)^\alpha}{\alpha} + \frac{u_2^2}{2},\qquad u=(u_1,u_2)\in U.
 \end{align*}
Indeed, in case of \ $\alpha\in(1,2)$, \ the formula for \ $R(u)$, $u\in U$, \ can be checked as follows.
By Dawson and Li \cite{DawLi}, in case of \ $\alpha\in(1,2)$,
 \begin{align*}
  R(u) &:= -bu_1 +\frac{u_2^2}{2}
         + \int_0^\infty \int_{-\infty}^\infty (\ee^{\langle u,\xi\rangle} - 1 - \langle u,\xi\rangle)
           C_\alpha\xi_1^{-1-\alpha}\,\dd\xi_1\delta_0(\dd\xi_2)  \\
       & = -bu_1 +\frac{u_2^2}{2}
         + C_\alpha \int_0^\infty (\ee^{u_1\xi_1} - 1 - u_1\xi_1)\xi_1^{-1-\alpha}\,\dd\xi_1\\
       & = -bu_1  +\frac{(-u_1)^\alpha}{\alpha} +\frac{u_2^2}{2},
         \qquad u\in U,
 \end{align*}
 where for the last equality we used that \ $1/\Gamma(-\alpha) = \alpha(\alpha-1)/\Gamma(2-\alpha)$, \ the imaginary part
 of \ $-i u_1\xi_1$ \ is non-negative and that
 \[
   (-iz)^\alpha = \frac{\alpha(\alpha-1)}{\Gamma(2-\alpha)}\int_0^\infty (\ee^{iz\xi_1} - 1 - iz\xi_1) \xi_1^{-1-\alpha}\,\dd\xi_1,
 \]
 for all complex numbers \ $z$ \ with non-negative imaginary part, see, e.g., Zolotarev \cite[pages 67 and 68]{Zol}.

Note also that for all \ $u=(u_1,u_2)\in U$ \ and \ $t\geq 0$, \ the real part of \ $\psi^{(1)}_t(u)$ \ is less than or equal to \ $0$
 \ (compare also with Remark 2.2 in Duffie et al. \cite{DufFilSch}), and, in addition, if \ $u_1\in\RR$ \ such that \ $u_1\leq 0$, \
 then \ $\psi^{(1)}_t(u)\in\RR$ \ with \ $\psi^{(1)}_t(u)\leq 0$.
\ Moreover, for all \ $t\geq 0$, \ we have
 \begin{align*}
  \phi_t(u) = \int_0^t (a\psi_s^{(1)}(u) + m\ee^{-\theta s}u_2)\,\dd s
            = a \int_0^t \psi_s^{(1)}(u)\,\dd s + mu_2\frac{1-\ee^{-\theta t}}{\theta}.
 \end{align*}

In fact, one can give a simplified characterization of the transition semigroup \ $(P_t)_{t\geq 0}$ \ by
 \begin{align}\label{help15}
  \int_0^\infty\int_{-\infty}^\infty & \ee^{-\lambda_1\xi_1 + i\lambda_2\xi_2} P_t((y_0,x_0),\dd\xi_1,\dd\xi_2)
     = \exp\left\{ -y_0 v_t(\lambda_1,\lambda_2) + ix_0 \ee^{-\theta t}\lambda_2 + g_t(\lambda_1,\lambda_2)\right\}
 \end{align}
 for \ $(\lambda_1,\lambda_2), \, (y_0,x_0) \in\RR_+\times\RR$, \ where
 \begin{align}\label{help53}
    g_t(\lambda_1,\lambda_2)
                 := \int_0^t( -av_s(\lambda_1,\lambda_2) + im\ee^{-\theta s}\lambda_2)\,\dd s
                  = -a \int_0^t v_s(\lambda_1,\lambda_2)\,\dd s + im\lambda_2 \frac{1-\ee^{-\theta t}}{\theta},
 \end{align}
 and \ $v_t(\lambda_1,\lambda_2)$, $t\geq 0$, \ is the unique non-negative
 solution of the differential equation
 \begin{align}\label{help13}
  \begin{cases}
   \frac{\partial v_t}{\partial t} (\lambda_1,\lambda_2)
       = -b v_t(\lambda_1,\lambda_2) - \frac{1}{\alpha}(v_t(\lambda_1,\lambda_2))^\alpha
          + \frac{1}{2}\ee^{-2\theta t}\lambda_2^2,\qquad t\geq 0,\\
   v_0(\lambda_1,\lambda_2) = \lambda_1
  \end{cases}
 \end{align}
 in case \ $\alpha\in(1,2]$.
\ Indeed, by \eqref{help12} with \ $u_1:=-\lambda_1$ \ and \ $u_2:=i\lambda_2$, \ we have
 \begin{align*}
  \int_0^\infty&\int_{-\infty}^\infty \ee^{-\lambda_1\xi_1 + i\lambda_2\xi_2}P_t((y_0,x_0),\dd\xi_1,\dd\xi_2)\\
    & = \exp\left\{ y_0\psi^{(1)}_t(-\lambda_1,i\lambda_2) + i x_0 \ee^{-\theta t}\lambda_2
                    + a \int_0^t \psi^{(1)}_s(-\lambda_1,i\lambda_2)\,\dd s
                    + im\lambda_2 \frac{1-\ee^{-\theta t}}{\theta} \right\},
 \end{align*}
 where
 \begin{align*}
  \frac{\partial \psi^{(1)}_t}{\partial t} (-\lambda_1,i\lambda_2)
      &= -b \psi^{(1)}_t(-\lambda_1,i\lambda_2)
         + \frac{1}{\alpha}(-\psi^{(1)}_t(-\lambda_1,i\lambda_2))^\alpha
         + \frac{1}{2}\ee^{-2\theta t}(i\lambda_2)^2,\quad t\geq 0,
 \end{align*}
 with \ $\psi^{(1)}_0(-\lambda_1,i\lambda_2) = -\lambda_1$ \ in case \ $\alpha\in(1,2]$.
\ Introducing the notation \ $v_t(\lambda_1,\lambda_2):=-\psi^{(1)}_t(-\lambda_1,i\lambda_2)$, \ we have
 the differential equation \eqref{help13} for \ $\alpha\in(1,2].$
\ Note also that
 \[
     v_t(\lambda_1,\lambda_2)\geq 0\qquad  \text{for all \ $(\lambda_1,\lambda_2)\in\RR_+\times\RR$,}
 \]
 since \ $\psi^{(1)}_t(u)\leq 0$ \ for \ $u_1\leq 0$.
\ The uniqueness of the solutions of the differential equation (Cauchy problem) \eqref{help13}
 follows by general results of Duffie et al. \cite[Propositions 6.1, 6.4 and Lemma 9.2]{DufFilSch}.
As follows, we present another direct proof based on the global version of the Picard--Lindel\"of existence and uniqueness theorem:
 if \ $f:D\to\RR$ \ is a continuous function on a connected, open set \ $D\subseteq \RR^2$ \ satisfying the global Lipschitz
 condition in its second variable
 \[
   \vert f(t,z_1) - f(t,z_2)\vert \leq C\vert z_1 - z_2\vert, \qquad (t,z_1),(t,z_2)\in D,
 \]
 with some constant \ $C>0$, \ then for all \ $(t_0,z_0)\in D$, \ the Cauchy problem
 \[
    z'(t) = f(t,z(t)) \qquad \text{with initial value \ $z(t_0) = z_0$}
 \]
 has a unique solution \ $z$ \ defined on a maximal interval of the form \ $(t_{-},t_{+})$, \ where \ $t_{-}<t_{+}$,
 \ $t_{-}\in[-\infty,\infty)$, $t_{+}\in(-\infty,\infty]$, \ and \ $(t,z(t))$ \ leaves any compact subset of \ $D$ \
 as \ $t\downarrow t_{-}$ \ and \ $t\uparrow t_{+}$, \ see, e.g., Hartman \cite[Chapter II, Theorems 1.1 and 3.1]{Har}
  and Walter \cite[Chapter 3]{Wal}.
Given \ $(\lambda_1,\lambda_2)\in\RR_+\times\RR$, \ for all \ $n\geq\lambda_1$, $n\in\NN$, \ let
 \ $D_n:=(-1,\infty)\times (-n,n)$ \ and let \ $f_n:D_n\to\RR$ \ be defined by
 \begin{align*}
  f_n(t,z)
   := -bz - \frac{1}{\alpha}z^\alpha + \frac{1}{2}\ee^{-2\theta t}\lambda_2^2, \qquad (t,z)\in D_n
   \qquad \text{ with \ $\alpha\in(1,2]$.}
 \end{align*}
Further, let \ $t_0:=0$ \ and \ $z_0:=\lambda_1$.
\ Then \ $D_n$ \ is open and connected, \ $f_n$ \ is continuous and satisfies the global Lipschitz condition in its second variable
 since, by the mean value theorem, for all \ $(t_1,z_1), (t_2,z_2)\in D$ \ we get
 \begin{align*}
   \vert f_n(t_1,z_1) - f_n(t_2,z_2) \vert
     \leq  b\vert z_1 - z_2\vert + \frac{1}{\alpha}\vert z_1^\alpha - z_2^\alpha\vert
           \leq b\vert z_1 - z_2\vert + n^{\alpha-1}\vert z_1 - z_2\vert.
 \end{align*}
By the global version of the Picard--Lindel\"of existence and uniqueness theorem, there is a unique solution \ $z_n$ \ of the
 Cauchy problem
 \[
    z_n'(t) = f_n(t,z_n(t)) \qquad \text{with initial value \ $z_n(t_0) = z_0$,}
 \]
 defined on a maximal interval of the form \ $((t_n)_{-},(t_n)_{+})$, \ where \ $(t_n)_{-} < (t_n)_{+}$, \ $(t_n)_{-}\in[-\infty,\infty)$ \ and
 \ $(t_n)_{+}\in(-\infty,+\infty]$.
\ Further, the solution \ $z_n$ \ leaves every compact subset of \ $D_n$ \ which implies that \ $(t_n)_{+}=\infty$ \ for all \ $n\in\NN$.
\ This yields the uniqueness of the solutions of the differential equation (Cauchy problem) \eqref{help13}.

{\bf Step 2.}
We show that
 \begin{align}\label{help17}
   v_t(\lambda_1,\lambda_2) \leq M(\lambda_1,\lambda_2) (1+t)\max(\ee^{-2\theta t},\ee^{-b t}),
   \qquad t\geq 0,\; (\lambda_1,\lambda_2)\in\RR_+\times\RR,
 \end{align}
 where
 \[
    M(\lambda_1,\lambda_2):=\begin{cases}
                              \lambda_1 + \frac{\lambda_2^2}{2\vert b-2\theta\vert}
                                  & \text{if $b\ne 2\theta$,}\\[2mm]
                              \lambda_1 + \frac{\lambda_2^2}{2}
                                     & \text{if $b = 2\theta$.}
                            \end{cases}
 \]
The proof is based on the following version of comparison theorem (see, e.g., Volkmann \cite{Vol}
 or Lemma B.3. in Filipovi\'{c} et al. \cite{FilMaySch}): if \ $S:\RR_+\times\RR\to\RR$ \ is a continuous function
 which is locally Lipschitz continuous in its second variable and \ $p,q:\RR_+\to\RR$ \ are differentiable functions
 satisfying
 \begin{align*}
   & p'(t) \leq S(t,p(t)),\quad t\geq 0,\\
   & q'(t) = S(t,q(t)),\quad t\geq 0,\\
   & p(0)\leq q(0),
 \end{align*}
 then \ $p(t)\leq q(t)$ \ for all \ $t\geq 0$.
\ Using this one can check that \ $v_t(\lambda_1,\lambda_2)\leq u_t(\lambda_1,\lambda_2)$ \ for all
 \ $t\geq 0$ \ and \ $(\lambda_1,\lambda_2)\in\RR_+\times\RR$, \ where for all \ $(\lambda_1,\lambda_2)\in\RR_+\times\RR$,
 \ $u_t(\lambda_1,\lambda_2)$, $t\geq 0$, \ is the unique solution of the differential equation
 \begin{align}\label{help16}
  \begin{cases}
   \frac{\partial u_t}{\partial t} (\lambda_1,\lambda_2)
       = -b u_t(\lambda_1,\lambda_2)
          + \frac{1}{2}\ee^{-2\theta t}\lambda_2^2,\qquad t\geq 0,\\
   u_0(\lambda_1,\lambda_2) = \lambda_1.
  \end{cases}
 \end{align}
Further, one can verify that
 \begin{align}\label{help22}
 u_t(\lambda_1,\lambda_2) =
  \begin{cases}
   \left( \lambda_1 + \frac{\lambda_2^2}{2(-b+2\theta)}\right)\ee^{-b t}
       - \frac{\lambda_2^2}{2(-b+2\theta)}\ee^{-2\theta t}
     & \text{if \ $b\ne 2\theta$,}\\[3mm]
     \left(\lambda_1  + \frac{\lambda_2^2}{2}t\right)\ee^{-b t}
    & \text{if \ $b = 2\theta$.}
  \end{cases}
 \end{align}
Indeed, the general solution of the homogeneous differential equation
 \ $\frac{\partial \widetilde u_t}{\partial t} (\lambda_1,\lambda_2) = -b \widetilde u_t(\lambda_1,\lambda_2)$, $t\geq 0$, \
 takes the form \ $\widetilde u_t(\lambda_1,\lambda_2) = C\ee^{-bt}$, $t\geq 0$, \ where \ $C\in\RR$, \ and
 it can be checked that a particular solution of the inhomogeneous differential equation \eqref{help16}
 (without the initial condition) is
 \[
   u_t(\lambda_1,\lambda_2)
     = \begin{cases}
         -\frac{\lambda_2^2}{2 (-b + 2\theta)}\ee^{-2\theta t} & \text{if \ $b\ne 2\theta$,}\\[1mm]
         \frac{\lambda_2^2}{2}t\ee^{-2\theta t} & \text{if \ $b=2\theta$.}
       \end{cases}
 \]
Hence the general solution of the differential equation \eqref{help16} (without the initial condition) takes the form
 \[
   u_t(\lambda_1,\lambda_2)
     = C\ee^{-bt}
      +\begin{cases}
         -\frac{\lambda_2^2}{2 (-b + 2\theta)}\ee^{-2\theta t} & \text{if \ $b\ne2\theta$,}\\[1mm]
         \frac{\lambda_2^2}{2}t\ee^{-2\theta t} & \text{if \ $b=2\theta$.}
      \end{cases}
 \]
Then taking into account the initial condition \ $u_0(\lambda_1,\lambda_2)=\lambda_1$, \ we have \eqref{help22}.

Finally, by \eqref{help22}, we readily have \eqref{help17} in case of \ $b=2\theta$.
If \ $b>2\theta$, \ then
 \begin{align*}
   u_t(\lambda_1,\lambda_2)
     \leq \lambda_1\ee^{-2\theta t}
           - \frac{\lambda_2^2}{2(-b+2\theta)}\ee^{-2\theta t}
     \leq \left( \lambda_1 - \frac{\lambda_2^2}{2(-b+2\theta)} \right) (1+t)\ee^{-2\theta t},
 \end{align*}
 and if \ $0<b<2\theta$, \ then
 \begin{align*}
   u_t(\lambda_1,\lambda_2)
     \leq \left( \lambda_1 + \frac{\lambda_2^2}{2(-b+2\theta)} \right) \ee^{-bt}
     \leq \left( \lambda_1 + \frac{\lambda_2^2}{2(-b+2\theta)} \right)(1+t) \ee^{-bt},
 \end{align*}
 as desired.

{\bf Step 3.}
By the continuity theorem and \eqref{help15}, to prove (i), it is enough to check that
 for all \ $(\lambda_1,\lambda_2), \,(y_0,x_0) \in\RR_+\times\RR$, \
 \begin{align}\label{help63}
  \begin{split}
   \lim_{t\to\infty} ( -y_0 v_t(\lambda_1,\lambda_2) + ix_0 \ee^{-\theta t}\lambda_2 + g_t(\lambda_1,\lambda_2) )
       &= -a\int_0^\infty v_s(\lambda_1,\lambda_2)\,\dd s + i\frac{m}{\theta}\lambda_2\\
       &=:g_\infty(\lambda_1,\lambda_2),
  \end{split}
 \end{align}
 and that the function \ $\RR_+\times\RR\ni(\lambda_1,\lambda_2) \mapsto g_\infty(\lambda_1,\lambda_2)$ \ is continuous.
Indeed, using \eqref{help15} and the independence of \ $(\eta_0,\zeta_0)$ \ and \ $(L_t,B_t)_{t\geq 0}$, \
 the law of total expectation yields that
 \begin{align*}
  \EE&(\ee^{-\lambda_1 Y_t + i\lambda_2 X_t})
      = \int_0^\infty \int_{-\infty}^\infty
       \EE\Big(\ee^{-\lambda_1 Y_t + i\lambda_2 X_t}\;\big\vert\; (Y_0,X_0) = (y_0,x_0) \Big)
        \,\PP_{(Y_0,X_0)}(\dd y_0,\dd x_0)\\
    & = \int_0^\infty \int_{-\infty}^\infty
       \exp\left\{-y_0 v_t(\lambda_1,\lambda_2) + ix_0\ee^{-\theta t}\lambda_2 + g_t(\lambda_1,\lambda_2)\right\}
               \,\PP_{(Y_0,X_0)}(\dd y_0,\dd x_0)
 \end{align*}
 for all \ $ (\lambda_1,\lambda_2)\in\RR_+\times\RR$, \ where \ $\PP_{(Y_0,X_0)}$ \ denotes the distribution of
 \ $(Y_0,X_0)$ \ on \ $\RR_+\times\RR$, \ and hence \eqref{help63} and the dominated convergence theorem implies that
 \begin{align*}
  \lim_{t\to\infty} \EE(\ee^{-\lambda_1 Y_t + i\lambda_2 X_t})
    =  \int_0^\infty \int_{-\infty}^\infty \ee^{g_\infty(\lambda_1,\lambda_2)}
                         \,\PP_{(Y_0,X_0)}(\dd y_0,\dd x_0)
    = \ee^{g_\infty(\lambda_1,\lambda_2)}
 \end{align*}
 for \ $(\lambda_1,\lambda_2)\in\RR_+\times\RR$.
\ Then, using the continuity of the function \ $\RR_+\times\RR\ni(\lambda_1,\lambda_2) \mapsto g_\infty(\lambda_1,\lambda_2)$ \
 (which will be checked later on), the continuity theorem yields (i).

Next we turn to prove \eqref{help63}.
Since \ $\theta>0$ \ and \ $b>0$, \ by \eqref{help17}, using also that \ $v_t(\lambda_1,\lambda_2)\geq 0$ \ for all \ $t\geq 0$
 \ and \ $(\lambda_1,\lambda_2)\in\RR_+\times\RR$ \ (see Step 1), we have
 \[
     \lim_{t\to\infty} (-y_0 v_t(\lambda_1,\lambda_2) + ix_0 \ee^{-\theta t}\lambda_2 ) = 0.
 \]
Recall that
 \begin{align*}
  g_t(\lambda_1,\lambda_2)
    = -a \int_0^t v_s(\lambda_1,\lambda_2)\,\dd s + im\lambda_2 \frac{1-\ee^{-\theta t}}{\theta}.
 \end{align*}
Since \ $\theta>0$, \ we have \ $\lim_{t\to\infty}\frac{1-\ee^{-\theta t}}{\theta} = \frac{1}{\theta}$, \
 and since \ $v_t(\lambda_1,\lambda_2)\geq 0$ \ for all \ $t\geq0$, $(\lambda_1,\lambda_2)\in\RR_+\times\RR$,
 \ (see Step 1), by dominated convergence theorem and \eqref{help17}, we get
 \[
   \lim_{t\to\infty} \int_0^t v_s(\lambda_1,\lambda_2)\,\dd s
      = \int_0^\infty v_s(\lambda_1,\lambda_2)\,\dd s.
 \]
Indeed, for all \ $t\geq 0$ \ and \ $s\geq 0$, \
 \begin{align*}
  \vert v_s(\lambda_1,\lambda_2)\mathbf 1_{[0,t]}(s)\vert
    \leq v_s(\lambda_1,\lambda_2),
 \end{align*}
 and, by \eqref{help17},
 \begin{align*}
    \int_0^\infty v_s(\lambda_1,\lambda_2)\,\dd s
      & \leq M(\lambda_1,\lambda_2)\int_0^\infty (1+s)\max(\ee^{-2\theta s},\ee^{-bs})\,\dd s\\
      & = M(\lambda_1,\lambda_2)\max\left(\frac{1}{b} + \frac{1}{b^2} , \frac{1}{2\theta} + \frac{1}{(2\theta)^2}\right) <\infty.
 \end{align*}

The continuity of the function \ $\RR_+\times\RR\ni(\lambda_1,\lambda_2) \mapsto g_\infty(\lambda_1,\lambda_2)$ \ can be checked as follows.
It will follow if we prove that for all \ $s\geq 0$, \ the function \ $v_s$ \ is continuous.
\ Namely, if \ $\lambda^{(n)}=(\lambda_1^{(n)},\lambda_2^{(n)})$, $n\in\NN$, \ is a sequence in \ $\RR_+\times\RR$, \ such that
 \ $\lim_{n\to\infty} \lambda^{(n)} = \lambda$, \ where \ $\lambda\in\RR_+\times\RR$, \ then
 \ $\lim_{n\to\infty}v_s(\lambda^{(n)})= v_s(\lambda)$ \ for all \ $s\geq 0$, \ and, by \eqref{help17},
 \[
    v_s(\lambda^{(n)}) = v_s(\lambda_1^{(n)},\lambda_2^{(n)})
                       \leq M(\lambda_1^{(n)},\lambda_2^{(n)})(1+s)\max(\ee^{-2\theta s},\ee^{-bs}),
                       \qquad n\in\NN,\;s\geq 0.
 \]
Since the sequence \ $\lambda^{(n)}$, $n\in\NN$, \ is bounded (being convergent), we have
 \[
   \sup_{n\in\NN} M(\lambda_1^{(n)},\lambda_2^{(n)}) <\infty,
 \]
 and using also that \ $\int_0^\infty (1+s)\max(\ee^{-2\theta s},\ee^{-bs})\,\dd s<\infty$, \ the dominated convergence theorem implies that
 \[
   \lim_{n\to\infty} \int_0^\infty v_s(\lambda_1^{(n)},\lambda_2^{(n)})\,\dd s
       = \int_0^\infty v_s(\lambda_1,\lambda_2)\,\dd s,
 \]
 which shows the continuity of \ $g_\infty$.
\ Finally, we turn to prove that for all \ $s\geq 0$, \ the function
 \ $\RR_+\times\RR\ni (\lambda_1,\lambda_2)\mapsto v_s(\lambda_1,\lambda_2)$ \ is continuous.
\ Note that the function \ $v_s$ \ does not depend on the parameters \ $a$ \ and \ $m$, \ since it is the unique solution of
 the differential equation \eqref{DE1}.
Let \ $(\tY_t, \tX_t)_{t\geq0}$ \ be an affine process satisfying the SDE
 \eqref{2dim_affine} with initial value \ $(\tY_0, \tX_0) = (Y_0,X_0)$ \ and
 with parameters \ $a=m=0$ \ and the given \ $b > 0$ \ and \ $\theta > 0$.
\ Then, by \eqref{help15},
  \begin{align}\label{help64}
  \int_0^\infty\int_{-\infty}^\infty \ee^{-\lambda_1\xi_1 + i\lambda_2\xi_2} \tP_s((y_0,x_0),\dd\xi_1,\dd\xi_2)
     = \exp\left\{ -y_0 v_s(\lambda_1,\lambda_2) + ix_0 \ee^{-\theta s}\lambda_2 \right\}
 \end{align}
 for \ $s\in\RR_+$, \ $(\lambda_1,\lambda_2), \, (y_0,x_0) \in\RR_+\times\RR$, \ where \ $(\tP_t)_{t\geq0}$ \ denotes the transition
 semigroup of the affine process \ $(\tY_t, \tX_t)_{t\geq0}$.
\ For all \ $s\in\RR_+$, \ the left-hand side of \eqref{help64} is continuous as a function of \ $(\lambda_1,\lambda_2)\in\RR_+\times\RR$, \ since
 for all \ $(\lambda_1,\lambda_2)\in\RR_+\times\RR$ \ and \ $(h_1,h_2)\in\RR^2$,
 \begin{align*}
  &\left\vert
    \int_0^\infty\!\!\int_{-\infty}^\infty \!\!\ee^{-(\lambda_1+h_1)\xi_1 + i(\lambda_2+h_2)\xi_2}
     \, \tP_s((y_0,x_0),\dd\xi_1,\dd\xi_2)
      - \int_0^\infty\!\!\int_{-\infty}^\infty \!\!\ee^{-\lambda_1\xi_1 + i\lambda_2\xi_2}
     \, \tP_s((y_0,x_0),\dd\xi_1,\dd\xi_2)
   \right\vert\\
  &\qquad \leq \int_0^\infty\int_{-\infty}^\infty \vert \ee^{-h_1\xi_1 + ih_2\xi_2} - 1\vert \, \tP_s((y_0,x_0),\dd\xi_1,\dd\xi_2) ,
 \end{align*}
 which tends to \ $0$ \ as \ $(h_1,h_2)\to (0,0)$, \ by dominated convergence theorem.
This implies that the right hand side of \eqref{help64} is also a continuous function \ $(\lambda_1,\lambda_2)\in\RR_+\times\RR$.
\ This readily yields the continuity of the function \ $v_s$ \ for all fixed \ $s\in\RR_+$.

{\bf (ii):} First we check that the one-dimensional distributions of \ $(Y_t,X_t)_{t\geq 0}$ \ are translation invariant and have common
 distribution as \ $(Y_\infty,X_\infty)$ \ has.
Using \eqref{help18}, \eqref{help15}, the tower rule and the independence of \ $(Y_0,X_0)$ \ and \ $(L,B)$, \ it is enough to check
 that for all \ $t\geq 0$ \ and \ $(\lambda_1,\lambda_2)\in\RR_+\times\RR$,
 \begin{align*}
   \EE\left(\exp\left\{ - v_t(\lambda_1,\lambda_2)Y_\infty + i \ee^{-\theta t}\lambda_2 X_\infty + g_t(\lambda_1,\lambda_2)\right\}\right)
     =\exp\left\{-a\int_0^\infty v_s(\lambda_1,\lambda_2)\,\dd s + i\frac{m}{\theta}\lambda_2\right\}.
 \end{align*}
By \eqref{help18}, \eqref{help53} and using also that \ $v_t(\lambda_1,\lambda_2)\geq 0$ \ for all \ $t\geq 0$ \ and \ $(\lambda_1,\lambda_2)\in\RR_+\times\RR$
 \ (see Step 1 of the proof of part (i)), we have
 \begin{align*}
           &\EE\left(\exp\left\{ - v_t(\lambda_1,\lambda_2) Y_\infty + i \ee^{-\theta t}\lambda_2 X_\infty + g_t(\lambda_1,\lambda_2)\right\}\right)\\
           &\qquad=\exp\left\{-a\int_0^\infty v_s( v_t(\lambda_1,\lambda_2) , \ee^{-\theta t}\lambda_2)\,\dd s + i\frac{m}{\theta}\ee^{-\theta t}\lambda_2
                        + g_t(\lambda_1,\lambda_2)\right\}\\
           &\qquad=\exp\left\{-a \left(\int_0^\infty v_s( v_t(\lambda_1,\lambda_2) , \ee^{-\theta t}\lambda_2)\,\dd s + \int_0^t v_s(\lambda_1,\lambda_2)\,\dd s\right)
                        + i\frac{m}{\theta}\lambda_2 \right\}.
 \end{align*}
Hence it remains to check that
 \[
    \int_0^\infty v_s(\lambda_1,\lambda_2)\,\dd s
       = \int_0^\infty v_s( v_t(\lambda_1,\lambda_2) , \ee^{-\theta t}\lambda_2)\,\dd s + \int_0^t v_s(\lambda_1,\lambda_2)\,\dd s
 \]
 for \ $t\geq 0$, $(\lambda_1,\lambda_2)\in\RR_+\times\RR$, \ i.e.,
 \[
    \int_t^\infty v_s(\lambda_1,\lambda_2)\,\dd s
       = \int_0^\infty v_s( v_t(\lambda_1,\lambda_2) , \ee^{-\theta t}\lambda_2)\,\dd s,
       \qquad t\geq 0,\;\; (\lambda_1,\lambda_2)\in\RR_+\times\RR.
 \]
For this it is enough to check that
 \[
    v_s( v_t(\lambda_1,\lambda_2) , \ee^{-\theta t}\lambda_2)
       = v_{s+t}(\lambda_1,\lambda_2),\qquad s,t\geq 0,\;\; (\lambda_1,\lambda_2)\in\RR_+\times\RR,
 \]
 or equivalently
 \begin{align}\label{help26}
    v_t( v_s(\lambda_1,\lambda_2) , \ee^{-\theta s}\lambda_2)
       = v_{s+t}(\lambda_1,\lambda_2),\qquad s,t\geq 0,\;\; (\lambda_1,\lambda_2)\in\RR_+\times\RR.
 \end{align}
By \eqref{DE1}, we have
 \begin{align*}
  \frac{\partial v_{s+t}}{\partial t} (\lambda_1,\lambda_2)
    = -bv_{s+t}(\lambda_1,\lambda_2) - \frac{1}{\alpha}(v_{s+t}(\lambda_1,\lambda_2))^\alpha
                  + \frac{1}{2}\ee^{-2\theta(s+t)}\lambda_2^2,
      \qquad t\geq 0,
 \end{align*}
 with initial condition \ $v_{s+0}(\lambda_1,\lambda_2) = v_s(\lambda_1,\lambda_2)$.
\ Note also that, again by \eqref{DE1},
 \begin{align*}
   \frac{\partial v_t}{\partial t}(v_s(\lambda_1,\lambda_2) , \ee^{-\theta s}\lambda_2)
      & = -bv_t(v_s(\lambda_1,\lambda_2) , \ee^{-\theta s}\lambda_2)
         - \frac{1}{\alpha}(v_t(v_s(\lambda_1,\lambda_2) , \ee^{-\theta s}\lambda_2))^\alpha\\
      &\phantom{=\,}
        + \frac{1}{2}\ee^{-2\theta t}(\ee^{-\theta s}\lambda_2)^2, \qquad t\geq 0,
 \end{align*}
 with initial condition \ $v_0(v_s(\lambda_1,\lambda_2) , \ee^{-\theta s}\lambda_2) = v_s(\lambda_1,\lambda_2)$.
\ Hence, for all \ $s \geq 0$, \ the left and right sides of \eqref{help26}, as functions of \ $t \geq 0$, \ satisfy
 the differential equation \eqref{DE1} with the initial value \ $v_s(\lambda_1,\lambda_2)$.
\ Since \eqref{DE1} has a unique solution for all non-negative initial values, we obtain \eqref{help26}.

Finally, the strict stationarity (translation invariance of the finite dimensional distributions) of
 \ $(Y_t,X_t)_{t\geq 0}$ \ follows by the chain's rule for conditional expectations using also that it is a time homogeneous Markov process.
\proofend

\section{Ergodicity}\label{Section_ergod}

Such as the existence of a unique stationary distribution, the question of ergodicity for an affine process
 is also in the focus of current investigations.

Recently, Sandri\'{c} \cite{San} has proved ergodicity of so called stable-like processes
 using the same technique that we applied.
Further, the ergodicity of the so-called \ $\alpha$-root process with \ $\alpha\in(1,2]$ \
 (see, the first SDE of \eqref{2dim_affine}) and some statistical applications were given in Li and Ma \cite{LiMa}.

The following result states the ergodicity of the affine diffusion process given by
 the SDE \eqref{2dim_affine} with \ $\alpha=2$.

\begin{Thm}\label{Thm_ergodic2}
Let us consider the 2-dimensional affine diffusion model \eqref{2dim_affine} with \ $\alpha=2$,
 \ $a>0$, \ $b>0$, $m\in\RR$, \ $\theta>0$, \ and with a random initial value
 \ $(\eta_0,\zeta_0)$ \ independent of \ $(L_t,B_t)_{t\geq 0}$
 \ satisfying \ $\PP(\eta_0\geq 0) = 1$.
\ Then, for all Borel measurable functions \ $f:\RR^2\to\RR$ \ such that
 \ $\EE\vert f(Y_\infty,X_\infty)\vert<\infty$, \ we have
 \begin{equation}\label{help_ergodic}
   \PP\left( \lim_{T\to\infty} \frac{1}{T}\int_0^T f(Y_s,X_s)\,\dd s
   = \EE f(Y_\infty,X_\infty) \right)=1,
 \end{equation}
 where the distribution of \ $(Y_\infty,X_\infty)$ \ is given by \eqref{help18}
 and \eqref{DE1} with \ $\alpha=2$.
\end{Thm}

\noindent{\bf Proof.}
We use the notations of Meyn and Tweedie \cite{MeyTwe2}, \cite{MeyTwe3}.
Using Theorem 6.1 (so called Foster-Lyapunov criteria) in Meyn and Tweedie \cite{MeyTwe3}, it is enough to check that
 \begin{itemize}
   \item[(a)] $(Y_t,X_t)_{t\geq 0}$ \ is a right process (defined on page 38 in Sharpe \cite{Sha});
   \item[(b)]
    all compact sets are petite
     for some skeleton chain (skeleton chains and petite sets are defined on pages 491, 500
    in Meyn and Tweedie \cite{MeyTwe2}, and page 550 in Meyn and Tweedie \cite{MeyTwe1}, respectively);
   \item[(c)]
    there exist \ $c,d\in\RR$ \ with \ $c>0$ \ such that the inequality
     \[
       (\cA_n V)(y, x) \leq - c V(y, x) + d , \qquad (y, x) \in O_n
     \]
    holds for all \ $n \in \NN$, \ where
    \ $O_n := \{(y, x) \in \RR_+ \times \RR : \|(y, x)\| < n\}$ \ for each
    \ $n \in \NN$,
    \begin{align}\label{help54}
     V(y, x) := (y - c_1)^2 + (x - c_2)^2,\qquad (y, x) \in \RR_+ \times \RR,
    \end{align}
    with some appropriate \ $c_1, c_2 \in \RR$,
    \ and \ $\cA_n$ \ denotes the extended generator of the process
    \ $(Y_t^{(n)}, X_t^{(n)})_{t\geq 0}$ \ given by
     \[
       (Y_t^{(n)}, X_t^{(n)})
       := \begin{cases}
           (Y_t, X_t), & \text{for \ $t < T_n$,} \\
           (0,n), & \text{for \ $t \geq T_n$,}
          \end{cases}
     \]
    where the stopping time \ $T_n$ \ is defined by
    \ $T_n := \inf\{ t \in \RR_+
                     : (Y_t, X_t) \in (\RR_+ \times \RR) \setminus O_n\}$.
   (Here we note that instead of \ $(0,n)$ \ we could have chosen any fixed state in \ $(\RR_+ \times \RR) \setminus O_n$,
    \ and we could also have defined \ $(Y_t^{(n)}, X_t^{(n)})_{t\geq 0}$ \ as the stopped process
     \ $(Y_{t\wedge T_n}, X_{t\wedge T_n})_{t\geq 0}$, \
    see Meyn and Tweedie \cite[page 521]{MeyTwe3}.)
 \end{itemize}
Indeed, then Theorem 6.1 in Meyn and Tweedie \cite{MeyTwe3} yields the
 exponential ergodicity of the process \ $(Y_t, X_t)_{t\in\RR_+}$, \ namely,
 there exist \ $\beta>0$ \ and \ $B \in \RR_+$ \ such that
 \[
   \sup_{|g|\leq V+1}
    \Big|\EE\big(g(Y_t,X_t) \mid (Y_0, X_0) = (y_0, x_0)\big) - \EE(g(Y_\infty, X_\infty))\Big|
   \leq B (V(y_0,x_0) + 1) \ee^{-\beta t}
 \]
 for all \ $t \in \RR_+$ \ and \ $(y_0, x_0) \in \RR_+ \times \RR$, \ where the supremum is running for
 Borel measurable functions \ $g:\RR_+\times\RR\to\RR$.
\ According to the discussion after Proposition 2.5 in Bhattacharya \cite{Bha}, this implies \eqref{help_ergodic}.
Here we also point out that, due to Bhattacharya \cite{Bha}, we do not have to assume that
 \ $\PP(Y_0>0)=1$ \ in order to prove \eqref{help_ergodic}.

To prove (a), it is enough to show that the process \ $(Y_t, X_t)_{t\in\RR_+}$
 \ is a (weak) Feller (see Meyn and Tweedy \cite[Section 3.1]{MeyTwe2}), strong Markov process with continuous sample paths,
 see, e.g., Meyn and Tweedy \cite[page 498]{MeyTwe2}.
According to Proposition 8.2 (or Theorem 2.7) in Duffie et al.~\cite{DufFilSch}, the process
 \ $(Y_t, X_t)_{t\geq 0}$ \ is a Feller Markov process.
Since \ $(Y_t, X_t)_{t\geq 0}$ \ has continuous sample paths almost surely (especially, it is c\`{a}dl\`{a}g),
 it is automatically a strong Markov process, see, e.g., Theorem 1 on page 56 in Chung \cite{Chu}.

To prove (b), in view of Proposition 6.2.8 in Meyn and Tweedy \cite{MeyTwe4},
 it is sufficient to show that the skeleton chain \ $(Y_n, X_n)_{n\in\ZZ_+}$ \ is
 irreducible with respect to the Lebesgue measure on \ $\RR_+ \times \RR$
 \ (see, e.g., Meyn and Tweedy \cite[page 520]{MeyTwe3}),
 and admits the Feller property.
The skeleton chain \ $(Y_n, X_n)_{n\in\ZZ_+}$ \ admits the Feller property, since
 the process \ $(Y_t, X_t)_{t\geq 0}$ \ is a Feller process.
In order to check irreducibility of the skeleton chain \ $(Y_n, X_n)_{n\in\ZZ_+}$
 \ with respect to the Lebesgue measure on \ $\RR_+ \times \RR$,
 it is enough to prove that the conditional distribution of \ $(Y_1, X_1)$ \ given
 \ $(Y_0, X_0)$ \ is absolutely continuous (with respect to the Lebesgue
 measure on \ $\RR_+ \times \RR$) \ with a conditional density function
 \ $f_{{(Y_1, X_1)\mid (Y_0,X_0)}} : \RR^2 \times \RR^2 \to \RR_+$ \ such that
 \ $f_{{(Y_1, X_1)\mid (Y_0,X_0)}}(y, x \mid y_0, x_0) > 0$ \ for all
 \ $(y, x , y_0, x_0) \in (0, \infty) \times \RR \times \RR_+ \times \RR$.
\ Indeed, the Lebesgue measure on \ $\RR_+ \times \RR$ \ is \ $\sigma$-finite,
 and if \ $B$ \ is a Borel set in \ $\RR_+ \times \RR$ \ with positive Lebesgue
 measure, then
 \begin{align*}
  \EE\biggl( \sum_{n=0}^\infty \bone_B(Y_n, X_n)
             \,\bigg|\, (Y_0, X_0) = (y_0, x_0) \biggr)
  &\geq \PP((Y_1, X_1) \in B \mid (Y_0, X_0) = (y_0, x_0)) \\
  &= \iint\limits_B
      f_{{(Y_1, X_1)\mid (Y_0,X_0)}}(y, x \mid y_0, x_0) \, \dd y \, \dd x
    > 0
 \end{align*}
 for all \ $(y_0, x_0) \in \RR_+ \times \RR$.
\ The existence of \ $f_{{(Y_1, X_1)\mid (Y_0,X_0)}}$ \ with the required property can be checked as follows.
By Theorem \ref{Pro_affine}, we have
 \begin{align*}
   &Y_1 = \ee^{-b}\left( y_0 + a\int_0^1 \ee^{bu}\,\dd u
                             + \int_0^1 \ee^{bu}\sqrt{Y_u} \,\dd L_u\right),\\
   &X_1 = \ee^{-\theta}\left( x_0 + m\int_0^1 \ee^{\theta u}\,\dd u
                                 + \int_0^1 \ee^{\theta u}\sqrt{Y_u} \,\dd B_u\right),
 \end{align*}
 provided that \ $(Y_0,X_0) = (y_0,x_0)$, \ $(y_0, x_0) \in \RR_+ \times \RR$.
\ Recall that a two-dimensional random vector \ $\zeta$ \ is absolutely continuous if and only if
 \ $V\zeta + v$ \ is absolutely continuous for all invertable matrices
 \ $V \in \RR^{2 \times 2}$ \ and for all vectors \ $v\in\RR^2$, and
  if the density function of \ $\zeta$ \ is positive on a set \ $S \subset \RR^2$, \
  then the density function of \ $V\zeta + v$ \ is positive on the set \ $V S + v$.
\ Hence it is enough to check that the random vector
 \begin{equation}\label{rv}
   \left( \int_0^1 \ee^{bu}\sqrt{Y_u} \,\dd L_u,  \int_0^1 \ee^{\theta u}\sqrt{Y_u} \,\dd B_u \right)
 \end{equation}
 is absolutely continuous with respect to the Lebesgue measure on \ $\RR^2$ \ having a
 density function being strictly positive on the set
 \[
    \left\{y\in\RR : y > -y_0 - a\int_0^1 \ee^{bu}\,\dd u \right\}\times \RR.
 \]
\ For all \ $x\in\RR$ \ and \ $y\leq -y_0 - a\int_0^1 \ee^{bu}\,\dd u$, \ we have
 \begin{align*}
  &\PP\left( \int_0^1 \ee^{bu}\sqrt{Y_u} \,\dd L_u < y,
           \int_0^1 \ee^{\theta u}\sqrt{Y_u} \,\dd B_u < x \right) \\
  &\qquad = \PP\biggl( \ee^b Y_1 - y_0 - a\int_0^1 \ee^{bu}\,\dd u < y,
                          \int_0^1 \ee^{\theta u}\sqrt{Y_u} \,\dd B_u < x \biggr)
          \leq \PP(Y_1 < 0)
          = 0,
 \end{align*}
 since \ $\PP(Y_1 \geq 0) = 1$.
\ Note that the conditional distribution of \ $\int_0^1 \ee^{\theta u}\sqrt{Y_u} \,\dd B_u$
 \ given \ $(Y_t)_{t\in[0,1]}$ \ is a normal distribution with mean zero and with variance
 \ $\int_0^1 \ee^{2\theta u} Y_u \,\dd u$ \ due to the fact that \ $(Y_t)_{t\in[0,1]}$ \ and
 \ $(B_t)_{t\geq 0}$ \ are independent.
Indeed, \ $(Y_t)_{t\geq 0}$ \ is adapted to the augmented filtration corresponding to
 \ $\eta_0$ \ and \ $(L_t)_{t\geq 0}$ \ (see, e.g., Karatzas and Shreve \cite[page 285]{KarShr}),
 and using the independence of the standard Wiener processes \ $L$ \ and \ $B$, \ and
 Problem 2.7.3 in Karatzas and Shreve \cite{KarShr}, one can argue that this augmented filtration
 is independent of the filtration generated by \ $B$.
\ Hence, using again the independence of the standard Wiener processes \ $L$ \ and \ $B$,
 we get for all \ $x\in\RR$ \ and \ $y > -y_0 - a\int_0^1 \ee^{bu}\,\dd u$, \
 \begin{align*}
  &\PP\left( \int_0^1 \ee^{bu}\sqrt{Y_u} \,\dd L_u < y,
           \int_0^1 \ee^{\theta u}\sqrt{Y_u} \,\dd B_u < x \right) \\
  &= \PP\biggl( \ee^b Y_1 - y_0 - a\int_0^1 \ee^{bu}\,\dd u < y,
                          \int_0^1 \ee^{\theta u}\sqrt{Y_u} \,\dd B_u < x \biggr)\\
  &=\EE\left( \PP\left( Y_1 < \ee^{-b}\left( y + y_0 + a\int_0^1 \ee^{bu}\,\dd u\right),
                          \int_0^1 \ee^{\theta u}\sqrt{Y_u} \,\dd B_u < x
                          \,\Big\vert\, (Y_t)_{t\in[0,1]}\right)\right)\\
  &=\EE\left( \EE\left( \bone_{\{ Y_1 < \ee^{-b}\left( y + y_0 + a\int_0^1 \ee^{bu}\,\dd u\right) \}}
                        \bone_{\{ \int_0^1 \ee^{\theta u}\sqrt{Y_u} \,\dd B_u <x\}}
                        \,\Big\vert\, (Y_t)_{t\in[0,1]} \right) \right)
 \end{align*}
 \begin{align*}
  &=\EE\left( \bone_{\{ Y_1 < \ee^{-b}\left( y + y_0 + a\int_0^1 \ee^{bu}\,\dd u\right) \}}
              \EE\left( \bone_{\{ \int_0^1 \ee^{\theta u}\sqrt{Y_u} \,\dd B_u <x\}}
                        \,\Big\vert\, (Y_t)_{t\in[0,1]} \right) \right)\\
  &=\EE\left( \bone_{\{ Y_1 < \ee^{-b}\left( y + y_0 + a\int_0^1 \ee^{bu}\,\dd u\right) \}}
              \int_{-\infty}^x p\left(w;\int_0^1 \ee^{2\theta u} Y_u\,\dd u \right)\dd w
       \right),
 \end{align*}
  where
 \[
   p(w;\sigma^2) := \frac{1}{\sqrt{2\pi\sigma^2}} \ee^{-w^2/(2\sigma^2)} , \qquad
   w \in\RR, \quad \sigma > 0 .
 \]
Here we call the attention that due to the assumption \ $a>0$,
 \[
    \PP\left( \int_0^1 \ee^{2\theta u} Y_u \,\dd u >0 \right) = 1.
 \]
Then, by the law of total expectation, for all \ $x\in\RR$ \ and \ $y > -y_0 - a\int_0^1 \ee^{bu}\,\dd u$,
 \begin{align*}
   &\PP\left( \int_0^1 \ee^{bu}\sqrt{Y_u} \,\dd L_u < y,
           \int_0^1 \ee^{\theta u}\sqrt{Y_u} \,\dd B_u < x \right)\\
   &\quad= \int_0^{\ee^{-b}\left( y + y_0 + a\int_0^1 \ee^{bu}\,\dd u\right)}
     \int_{-\infty}^x
      \EE\left( \frac{1}{\sqrt{2\pi \int_0^1 \ee^{2\theta u} Y_u\,\dd u }}
                 \exp\left\{-\frac{w^2}{2 \int_0^1 \ee^{2\theta u} Y_u\,\dd u }\right\} \,\Bigg\vert\,  Y_1=z \right)
          f_{Y_1}(z)\,\dd z \,\dd w,
 \end{align*}
 where \ $f_{Y_1}$ \ denotes the density function of \ $Y_1$ \ (given that \ $Y_0=y_0$, \ $y_0\in\RR_+$).
\ In case of \ $y_0\in(0,\infty)$, \ it is given by
 \begin{align*}
   f_{Y_1}(y)
      = \frac{2b\ee^{b(2a+1)}}{\ee^b - 1}
        \left(\frac{y}{y_0}\right)^{a-1/2}
        \exp\left\{ - \frac{2b(y_0 + \ee^by)}{\ee^b - 1} \right\}
        I_{2a-1}\left(\frac{2b\sqrt{y_0y}}{\sinh(b/2)}\right)
        \bone_{(0,\infty)}(y), \quad y\in\RR,
 \end{align*}
 where \ $I_{2a-1}$ \ denotes the modified Bessel of the first kind of order \ $2a-1$, \ i.e.,
 \[
    I_{2a-1}(x) = \sum_{m=0}^\infty \frac{1}{m!\Gamma(m+2a)}\left(\frac{x}{2}\right)^{2m+2a-1},
    \qquad x\in(0,\infty),
 \]
 see, e.g., Cox et al. \cite[Equation (18)]{CoxIngRos}, Jeanblanc et al. \cite[Proposition 6.3.2.1]{JeaYorChe}
 or Ben Alaya and Kebaier \cite[the proof of Proposition 2]{BenKeb2}.
While in case of \ $y_0=0$, \ the density function of \ $Y_1$ \ (given that \ $Y_0=0$) \ is given by
 \[
   f_{Y_1}(y)
     = \frac{1}{\Gamma(2a)}\left(\frac{2b}{1-\ee^{-b}}\right)^{2a} y^{2a-1}
       \exp\left\{ - \frac{2by}{1-\ee^{-b}}\right\}\bone_{(0,\infty)}(y), \quad y\in\RR,
 \]
 since, by Ikeda and Watanabe \cite[page 222]{IkeWat},
 \ $Y_1$ \ (given that \ $Y_0=0$) \ has a Gamma distribution with parameters \ $2a$ \ and \ $2b/(1-\ee^{-b})$.
\ Note that in both cases \ $f_{Y_1}(y) > 0$ \ for all \ $y\in(0,\infty)$.
\ Then, by a change of variable, we have for all \ $x\in\RR$ \ and \ $y > -y_0 - a\int_0^1 \ee^{bu}\,\dd u$,
 \begin{align*}
   &\PP\left( \int_0^1 \ee^{bu}\sqrt{Y_u} \,\dd L_u < y,
           \int_0^1 \ee^{\theta u}\sqrt{Y_u} \,\dd B_u < x \right)\\
   &\quad = \int_{ -y_0 - a\int_0^1 \ee^{bu}\,\dd u}^y \int_{-\infty}^x
            \EE\left( \frac{\exp\Big\{-\frac{w^2}{2 \int_0^1 \ee^{2\theta u} Y_u\,\dd u }\Big\}}
                           {\sqrt{2\pi \int_0^1 \ee^{2\theta u} Y_u\,\dd u }}
                 \,\Bigg\vert\,  Y_1 = \ee^{-b}\left( v + y_0 + a\int_0^1 \ee^{bu}\,\dd u\right) \right)\\
   &\phantom{\quad = \int_{ -y_0 - a\int_0^1 \ee^{bu}\,\dd u}^y \int_{-\infty}^x\;\;}
         \times f_{Y_1}\left( \ee^{-b}\left( v + y_0 + a\int_0^1 \ee^{bu}\,\dd u\right) \right) \ee^{-b} \,\dd v \,\dd w.
 \end{align*}

In what follows we will make use of the following simple observation:
 if \ $\xi$ \ and \ $\eta$ \ are random variables such that \ $\PP(\xi>0)=1$,
 \ $\EE(\xi) < \infty$, \ $\PP(\eta>0)=1$, \ and \ $\eta$ \ is absolutely continuous
 with a density function \ $f_\eta$ \ having the property \ $f_\eta(x)>0$ \ Lebesgue a.e. \ $x\in(0,\infty)$, \
 then \ $\EE(\xi\mid \eta=y)>0$ \  Lebesgue a.e. \ $y\in(0,\infty)$.
\ For completeness, we give a proof.
Since the distribution \ $\PP_\eta$ \ of \ $\eta$ \ on \ $(\RR_+,\cB(\RR_+))$ \ (where \ $\cB(\RR_+)$ \ denotes
 the Borel $\sigma$-algebra on \ $\RR_+$) \ is equivalent to the Lebesgue meaure on \ $(0,\infty)$,
 \ it is equivalent to check that \ $\EE(\xi\mid \eta=y)>0$ \ $\PP_\eta$-a.e. \ $y\in(0,\infty)$.
\ By the definition of conditional expectation, it is also equivalent to check that
 \ $\EE(\xi\mid \eta)(\omega)>0$ \ $\PP$-a.e. \ $\omega\in\Omega$.
\ On the contrary, let us suppose that \ $\PP(\EE(\xi\mid \eta) =0)>0$.
\ Since \ $\PP(\EE(\xi\mid \eta) \geq 0)=1$ \ and \ $\EE(\xi\mid \eta)$ \ is measurable with respect to the
 \ $\sigma$-algebra \ $\sigma(\eta)$ \ generated by \ $\eta$, \ there would exist a set \ $A\in\sigma(\eta)$ \
 such that \ $\PP(A)>0$ \ and \ $\EE(\xi\mid \eta)(\omega)=0$ \ for all \ $\omega\in A$.
\ Then, using again the definition of conditional expectation, we should have
 \ $0=\EE( \EE(\xi\mid \eta) \bone_A) = \EE(\xi\bone_A)$.
\ Since \ $\PP(\xi\bone_A\geq 0)=1$, \ we would get \ $\PP(\{ \omega\in A : \xi(\omega)=0\})=1$,
 \ which leads us to a contradiction (due to that \ $\PP(\xi>0)=1$).

Now we turn back to the proof that the random vector \eqref{rv} is absolutely continuous
 with respect to the Lebesgue measure on \ $\RR^2$ \ with a density function being strictly positive on the set
 \ $\left\{y\in\RR : y > -y_0 - a\int_0^1 \ee^{bu}\,\dd u \right\}\times \RR$.
\ Since \ $\sup_{x\in\RR_+} x \ee^{-w^2x^2/2} < \infty$, $w\in\RR\setminus\{0\}$, \ implies
 \[
   \EE\left( \left(2\pi \int_0^1 \ee^{2\theta u} Y_u\,\dd u\right)^{1/2}
             \exp\Big\{-\frac{w^2}{2 \int_0^1 \ee^{2\theta u} Y_u\,\dd u }\Big\}
      \right)
   < \infty ,\quad w\in\RR\setminus\{0\},
 \]
 and using also that
 \ $f_{Y_1}(\ee^{-b} (v+y_0 + a\int_0^1 \ee^{bu}\,\dd u) ) >0$ \ for all \ $v>-y_0 - a\int_0^1 \ee^{bu}\,\dd u$, \
 there exists a measurable function \ $g:\RR^2\to\RR_+$ \ such that \ $g(v,w)>0$ \ for
 \ $v>-y_0 - a\int_0^1 \ee^{bu}\,\dd u$, \ $w\in\RR$, \ and
 \begin{align*}
  &\PP\left( \int_0^1 \ee^{bu}\sqrt{Y_u} \,\dd L_u < y,
           \int_0^1 \ee^{\theta u}\sqrt{Y_u} \,\dd B_u < x \right)\\
  &\quad=\begin{cases}
      \int_{ - y_0 - a\int_0^1 \ee^{bu}\,\dd u}^y \int_{-\infty}^x g(v,w)\,\dd v \,\dd w
            & \text{if \ $y > - y_0 - a\int_0^1 \ee^{bu}\,\dd u $, \ $x\in\RR$,}\\
      0 & \text{if \ $y \leq - y_0 - a\int_0^1 \ee^{bu}\,\dd u $, \ $x\in\RR$,}
    \end{cases}
 \end{align*}
 as desired.
Consequently, the random vector \eqref{rv} is absolutely continuous with density function \ $g$ \ having the desired property.

To prove (c), first we note that, since the sample paths of \ $(Y,X)$ \ are almost surely continuous,
 for each \ $n \in \NN$, \ the extended generator has the form
 \[
   (\cA_n f)(y,x)
   = \frac{1}{2} y f''_{1,1}(y,x) + \frac{1}{2} y f''_{2,2}(y,x)
     +(a - by) f'_1(y,x) + (m - \theta x) f'_1(y,x)
 \]
 for all \ $(y, x) \in O_n$ \ and \ $f \in \cC^2(\RR_+ \times \RR, \RR)$, \ see, e.g., page 538 in Meyn and Tweedie \cite{MeyTwe3}.
We also note that, by Duffie et al.~\cite[Theorem 2.7]{DufFilSch}, for functions
 \ $f \in \cC^2_{\mathrm{c}}(\RR_+ \times \RR, \RR)$, \ $\cA_n f = \cA f$ \ on \ $O_n$, \ where
 \ $\cA$ \ denotes the (non-extended) generator of the process
 \ $(Y_t, X_t)_{t\in\RR_+}$.
\ For the function \ $V$ \ defined in \eqref{help54},  we have \ $V\in\cC^2(\RR_+\times\RR, \RR)$ \ and
 \begin{align*}
   V_1'(y,x) = 2(y-c_1),\quad V_2'(y,x) = 2(x-c_1),
   \quad V_{1,1}''(y,x) = V_{2,2}''(y,x) = 2
 \end{align*}
 for \ $(y,x)\in\RR_+\times\RR$, \ and hence for all \ $n\in\NN$ \ and \ $0<c<\infty$,
 \begin{align*}
  &(\cA_nV)(y,x)  + c V(y,x)\\
  &\; = 2y + (a-by)2(y-c_1)
       + (m-\theta x)2(x-c_2)
       + c(y-c_1)^2 + c(x-c_2)^2 \\
  & \; = (c-2b)y^2 + 2(1+a+bc_1-cc_1)y + cc_1^2 -2ac_1
      + 2(m-\theta x)(x-c_2)
      + c(x-c_2)^2
 \end{align*}
 for all \ $(y,x)\in O_n$.
\ Let us choose \ $c>0$ \ in such a way that \ $c\ne 2b$ \ and let \ $c_2:=\frac{m}{\theta}$.
\ Then
 \begin{align*}
  &(\cA_nV)(y,x) + c V(y,x)\\
  &\quad = (c-2b)y^2 + 2(1+a+bc_1-cc_1)y + cc_1^2 -2ac_1
     - 2\theta \left(x-\frac{m}{\theta}\right)^2
     + c \left(x-\frac{m}{\theta}\right)^2 \\
  &\quad = (c-2b)\left(y + \frac{1+a+bc_1-cc_1}{c-2b}\right)^2
     - \frac{(1+a+bc_1-cc_1)^2}{c-2b}
     + cc_1^2 - 2ac_1\\
  &\phantom{\quad =\,}
     + (c-2\theta)\left(x - \frac{m}{\theta}\right)^2,
     \qquad (y,x)\in O_n.
 \end{align*}
Hence,
 \[
    (\cA_nV)(y,x) \leq -cV(y,x) + d, \qquad (y,x)\in O_n,\;\;n\in\NN,
 \]
 with \ $c\in(0,2\min(b,\theta))$ \ and
 \[
    d:=-\frac{(1+a+bc_1-cc_1)^2}{c-2b} + cc_1^2 - 2ac_1.
 \]
Note that the above argument also shows that in the definition \eqref{help54} of the function \ $V$ \
 the constant \ $c_1$ \ can be an arbitrary real number.
\proofend

In the next theorem we collected several facts about the limiting random variable \ $(Y_\infty,X_\infty)$ \
 given by \eqref{help18} and \eqref{DE1} with \ $\alpha=2$.

\begin{Thm}\label{Thm_ergodic3}
The random variable \ $(Y_\infty,X_\infty)$ \ given by \eqref{help18} and \eqref{DE1} with \ $\alpha=2$ \ is absolutely continuous,
 the Laplace transform of \ $Y_\infty$ \ takes the form
 \begin{align}\label{help74}
  \EE(\ee^{-\lambda_1 Y_\infty}) = \left(1 + \frac{\lambda_1}{2b}\right)^{-2a},
           \qquad \lambda_1\in\RR_+,
 \end{align}
 yielding that \ $Y_\infty$ \ has Gamma distribution with parameters \ $2a$ \ and \ $2b$.
\ Further, all the (mixed) moments of \ $(Y_\infty,X_\infty)$ \ of any order
 are finite, i.e., \ we have \ $\EE(Y_\infty^n\vert X_\infty\vert^p)<\infty$
 \ for all \ $n,p\in\ZZ_+$, \ and especially,
 \begin{align*}
          &\EE(Y_\infty)=\frac{a}{b}, \qquad \EE(X_\infty)=\frac{m}{\theta}, \\
          &\EE(Y_\infty^2) = \frac{a(2a+1)}{2b^2},
         \qquad
         \EE(Y_\infty X_\infty) = \frac{ma}{\theta b},
         \qquad
         \EE(X_\infty^2) = \frac{a\theta+ 2bm^2}{2b\theta^2},\\
         &\EE(Y_\infty X_\infty^2)
            = \frac{a}{(b+2\theta)2b^2\theta^2}
               \big(\theta(ab+2a\theta+\theta) + 2m^2b(2\theta + b)\big).
 \end{align*}
\end{Thm}

\noindent{\bf Proof.}
First we show that the random variable \ $(Y_\infty, X_\infty)$ \ is absolutely
 continuous.
Let us consider the 2-dimensional affine diffusion model \eqref{2dim_affine} with \ $\alpha=2$,
 \ $a>0$, \ $b>0$, $m\in\RR$, \ $\theta>0$, \ and with a random initial value
 \ $(Y_0,X_0)$ \ independent of \ $(L_t,B_t)_{t\geq 0}$ \ having the same distribution as that of \ $(Y_\infty,X_\infty)$.
\ Then, by part (ii) of Theorem \ref{Thm_ergodic1}, the process \ $(Y_t,X_t)_{t\geq 0}$ \ is strictly stationary.
Hence it is enough to prove that \ $(Y_1,X_1)$ \ is absolutely continuous.
This can be done similarly to the proof of condition (b) in the proof of Theorem \ref{Thm_ergodic2}, we only sketch a proof.
Namely, by Theorem \ref{Pro_affine}, we have
 \begin{align*}
   &Y_1 = \ee^{-b}\left( Y_0 + a\int_0^1 \ee^{bu}\,\dd u
                             + \int_0^1 \ee^{bu}\sqrt{Y_u} \,\dd L_u\right),\\
   &X_1 = \ee^{-\theta}\left( X_0 + m\int_0^1 \ee^{\theta u}\,\dd u
                                 + \int_0^1 \ee^{\theta u}\sqrt{Y_u} \,\dd B_u\right),
 \end{align*}
Using that if \ $\xi$ \ and \ $\eta$ \ are independent two-dimensional random vectors such that one of them is absolutely continuous,
 then their sum \ $\xi+\eta$ \ is absolutely continuous (regardless of the other distribution),
 it is enough to check that the random vector
 \begin{equation*}
   \left( \int_0^1 \ee^{bu}\sqrt{Y_u} \,\dd L_u,  \int_0^1 \ee^{\theta u}\sqrt{Y_u} \,\dd B_u \right)
 \end{equation*}
 is absolutely continuous with respect to the Lebesgue measure on \ $\RR^2$.
\ Its proof goes along the very same lines as in the proof of condition (b) in the proof of Theorem \ref{Thm_ergodic2}.

The fact that \ $Y_\infty$ \ has Gamma distribution with parameters \ $2a$ \ and \ $2b$ \ follows by Cox et al.
 \cite[Equation (20)]{CoxIngRos}.
In what follows we give two other proofs.
By \eqref{help18},
 \[
   \EE(\ee^{-\lambda_1 Y_\infty})
       = \exp\left\{ -a \int_0^\infty v_s(\lambda_1,0)\,\dd s\right\},
       \qquad \lambda_1\in\RR_+,
 \]
 where \ $v_t(\lambda_1,0)$, $t\geq 0$, \ is the unique solution of the differential equation
 \begin{align}\label{DE3}
   \begin{cases}
          \frac{\partial v_t}{\partial t}(\lambda_1,0)
               = -bv_t(\lambda_1,0) - \frac{1}{2}(v_t(\lambda_1,0))^2,\qquad t\geq 0,\\
          v_0(\lambda_1,0) = \lambda_1.
   \end{cases}
 \end{align}
Note that the formula above for \ $\EE(\ee^{-\lambda_1 Y_\infty})$ \ is a special case of formula (3.31) in Li \cite{Li}.
The differential equation \eqref{DE3} is of Bernoulli type, and hence for \ $\lambda_1>0$, \
 with the transformation \ $u_t(\lambda_1,0):=v_t(\lambda_1,0)^{-1}$, $t\geq 0$, \ we have
 the inhomogeneous linear differential equation with constant coefficients:
 \begin{align}\label{DE4}
   \begin{cases}
          \frac{\partial u_t}{\partial t}(\lambda_1,0)
               = bu_t(\lambda_1,0) + \frac{1}{2},\qquad t\geq 0,\\
          u_0(\lambda_1,0) = \lambda_1^{-1}.
   \end{cases}
 \end{align}
One can check that
 \[
   v_t(\lambda_1,0)
     = u_t(\lambda_1,0)^{-1}
     = \left( \left( \frac{1}{\lambda_1} + \frac{1}{2b}\right)\ee^{bt} - \frac{1}{2b} \right)^{-1},
     \qquad \lambda_1>0, \;\; t\geq 0.
 \]
Hence,
 \begin{align*}
   \EE(\ee^{-\lambda_1 Y_\infty})
   & = \exp\left\{ -a \int_0^\infty \frac{1}{\left( \frac{1}{\lambda_1} + \frac{1}{2b}\right)\ee^{bs} - \frac{1}{2b}} \,\dd s \right\}
     = \exp\left\{ -a \int_1^\infty \frac{1}{z\left(\left( \frac{b}{\lambda_1} + \frac{1}{2}\right)z - \frac{1}{2}\right)} \,\dd z \right\}\\
   & = \exp\left\{ -a \int_1^\infty \left( -\frac{2}{z} + \frac{1+\frac{2b}{\lambda_1}}
         { \left(\frac{b}{\lambda_1} + \frac{1}{2}\right)z- \frac{1}{2}} \right) \,\dd z \right\} \\
   & = -2a\lim_{z\to\infty} \left( -\ln z + \ln\left( \left(\frac{b}{\lambda_1} + \frac{1}{2}\right)z- \frac{1}{2} \right) \right)
       + 2a\ln\left( \frac{b}{\lambda_1}\right) \\
   & = -2a\ln\left( 1 + \frac{\lambda_1}{2b}\right),
 \end{align*}
 which yields \eqref{help74}.

Next we give another argument for deriving \eqref{help74}.
Since in case of \ $\alpha=2$, \ the process \ $(Y_t)_{t\geq 0}$ \ is a continuous state branching
 process with branching mechanism \ $bz+z^2/2$, $z\geq 0$, \
 and with immigration mechanism \ $az$, $z\geq 0$, \ by the proof of Theorem 3.20 in Li \cite{Li},
 we have
 \begin{align*}
  \int_0^\infty v_s(\lambda_1,0)\,\dd s
   = \int_0^{\lambda_1} \frac{az}{bz + z^2/2}\,\dd z
   = 2a\ln\left(1+\frac{\lambda_1}{2b}\right), \qquad \lambda_1\in\RR_+,
 \end{align*}
 which yields \eqref{help74}.

Now we prove that all the mixed moments of \ $(Y_\infty, X_\infty)$ \ are finite.
Let us consider now the 2-dimensional affine diffusion model \eqref{2dim_affine} with \ $\alpha=2$,
 \ $a>0$, \ $b>0$, \ $m\in\RR$, \ $\theta>0$, \ and with a random initial value
 \ $(Y_0,X_0)$ \ independent of \ $(L_t,B_t)_{t\geq 0}$ \ such that all the mixed moments of
 \ $(Y_0, X_0)$ \ are finite and \ $\PP(Y_0>0)=1$.
\ We note that, due to Theorem \ref{Thm_ergodic1}, the distribution of \ $(Y_\infty,X_\infty)$ \ does not depend on
 the initial value of the model \eqref{2dim_affine}, so we can have such a choice.
First we show that
 \begin{align}\label{help65}
    \int_0^t\EE(Y_u^n X_u^{2p})\,\dd u<\infty\qquad  \text{ for all \ $t\geq 0$ \ and \ $n,p\in\ZZ_+$.}
 \end{align}
Due to the inequality between two power means,
 \ $(a+b+c)^{2p}\leq K(a^{2p} + b^{2p} + c^{2p})$, $a,b,c\in\RR$, \ with some constant \ $K>0$, \
 and hence, by \eqref{help2}, to prove \eqref{help65} it is enough to check that
 \begin{align}\label{help75}
  \int_0^t \EE\Bigg( Y_u^n \Bigg(
                \ee^{-2p\theta u} X_0^{2p}
                 + \left( m \int_0^u \ee^{-\theta(u-v)}\,\dd v \right)^{2p}
                  +  \left( \int_0^u \ee^{-\theta(u-v)}\sqrt{Y_u }\,\dd B_v \right)^{2p}
             \Bigg) \Bigg)\dd u
       <\infty
 \end{align}
 for all \ $t\geq 0$ \ and \ $n,p\in\ZZ_+$.
\ Since for all \ $u\in[0,t]$ \ the conditional distribution of \ $\int_0^u \ee^{-\theta(u-v)}\sqrt{Y_v}\,\dd B_v$ \
 with respect to the \ $\sigma$-algebra generated by \ $(Y_s)_{s\in[0,t]}$ \ is a normal distribution with mean \ $0$
 \ and with variance \ $\int_0^u \ee^{-2\theta(u-v)} Y_v \,\dd v$, \ for proving \eqref{help75} it is enough to check that
 for all \ $t\geq 0$ \ and \ $n,p\in\ZZ_+$,
 \begin{align*}
     \int_0^t \EE(\ee^{-2p\theta u}Y_u^n X_0^{2p})\,\dd u<\infty,\qquad  \int_0^t \EE(Y_u^n)\,\dd u<\infty,
 \end{align*}
 and
 \begin{align*}
  \int_0^t \EE\left( Y_u^n \left(\int_0^u \ee^{-2\theta(u-v)}Y_v \,\dd v \right)^{p} \right)\,\dd u < \infty.
 \end{align*}
Using that \ $\theta>0$, \ $\PP(Y_t\geq 0, \; t\in\RR_+)=1$, \ by Cauchy-Schwarz inequality, for all \ $u\in[0,t]$, \ we have
 \ $\EE(\ee^{-2p\theta u}Y_u^n X_0^{2p})
      \leq \left(\EE(Y_u^{2n})\right)^{1/2}
           \left(\EE(X_0^{4p})\right)^{1/2}$, \ and
 \begin{align*}
   &\EE\left( Y_u^n \left(\int_0^u \ee^{-2\theta(u-v)}Y_v \,\dd v \right)^{p} \right)
     \leq \Big( \EE\big(Y_u^{2n} \big)\Big)^{1/2}
       \left( \EE \left( \left(\int_0^u Y_v \,\dd v \right)^{2p} \right) \right)^{1/2} \\
   &\qquad= \Big( \EE(Y_u^{2n} )\Big)^{1/2}
        \left( \int_0^u\cdots\int_0^u
              \EE\big(Y_{v_1}\cdots Y_{v_{2p}} \big)\,\dd v_1\cdots\dd v_{2p}
         \right)^{1/2} \\
  &\qquad \leq \Big( \EE\big(Y_u^{2n} \big)\Big)^{1/2}
          \left( \int_0^u\cdots\int_0^u
              \left( \EE\big(Y_{v_1}^{2p} \big)\cdots \EE\big(Y_{v_{2p}}^{2p} \big) \right)^{1/(2p)}\dd v_1\cdots\dd v_{2p}
         \right)^{1/2},
 \end{align*}
 where the last inequality follows by the multivariate version of H\"older's inequality.
Since \ $\EE(X_0^{4p})<\infty$, \ this shows that in order to prove \eqref{help65}
 it is enough to check that for all \ $k\in\ZZ_+$ \ and \ $t\in\RR_+$,
 \[
   \sup_{u\in[0,t]} \EE(Y_u^k) <\infty.
 \]
According to Proposition 3 in Ben Alaya and Kebaier \cite{BenKeb2},
 \begin{align}\label{help72}
   \sup_{s\geq 0} \EE(Y_s^k\mid Y_0=y_0)<\infty \qquad \text{for all \ $y_0>0$ \ and \ $k\in\ZZ_+$.}
 \end{align}
Next we show that for all \ $s\geq 0$, $y_0>0$, \ and \ $k\in\ZZ_+$,
 \[
   \EE(Y_s^k\mid Y_0=y_0) \quad \text{is a polynomial of \ $y_0$ \ of degree \ $k$.}
 \]
By \eqref{help15}, for all \ $s\geq 0$, $\lambda_1\geq 0$, \ and \ $(y_0,x_0)\in\RR_+\times\RR$,
 \begin{align}\label{help73}
  \int_0^\infty\int_{-\infty}^\infty \ee^{-\lambda_1\xi_1} P_s((y_0,x_0),\dd\xi)
     = \exp\left\{ -y_0 v_s(\lambda_1,0) - a \int_0^s v_u(\lambda_1,0)\,\dd u\right\}.
 \end{align}
Since, by \eqref{help72}, \ $\EE(Y_s^k\mid Y_0=y_0)<\infty$ \ for all \ $y_0>0$ \ and \ $k\in\ZZ_+$, \ and,
 by \eqref{DE1}, the function \ $v_s(\lambda_1,0)$, $s\geq 0$, \ is continuously differentiable of infinitely many times,
 one can differentiate both sides of \eqref{help73} with respect to \ $\lambda_1$ \ for \ $k$ \ times.
Since \ $v_s(0,0)=0$ \ for all \ $s\geq 0$ \ (which is a consequence of the uniqueness of the solution of
 the differential equation \eqref{DE1}), we have
 \[
   \EE(Y_s^k\mid Y_0=y_0)
      = (-1)^k\frac{\dd^k}{\dd\lambda_1^k}
        \exp\left\{ -y_0 v_s(\lambda_1,0) - a\int_0^s v_u(\lambda_1,0)\,\dd u\right\}
        \Bigg\vert_{\lambda_1=0}
 \]
 is a polynomial of \ $y_0$ \ of degree \ $k$, \ and the coefficients of this polynomial are continuous functions
 of \ $s$.
\ Since a continuous function on a compact set is bounded, we get that \ $\sup_{u\in[0,t]} \EE(Y_u^k \mid Y_0 = y_0)$ \
 can be bounded above by a polynomial of \ $y_0$ \ having degree \ $k$, \ say \ $Q_k(y_0)$ \ (this polynomial depends also on \ $t$).
Hence, by the law of total expectation,
 \begin{align*}
   \sup_{u\in[0,t]} \EE(Y_u^k)
         = \sup_{u\in[0,t]} \int_0^\infty \EE(Y_u^k\mid Y_0=y_0) \PP_{Y_0}(\dd y_0)
            \leq \int_0^\infty Q_k(y_0) \PP_{Y_0}(\dd y_0)
        = \EE(Q_k(Y_0))<\infty,
 \end{align*}
 where the last inequality follows by the assumption that all the (mixed) moments of \ $(Y_0,X_0)$ \ of any order are finite.

We note that for proving the finiteness of \ $\int_0^t\EE(Y_u^n X_u^{2p})\,\dd u$ \ we could have used part (i) of Theorem 2.16
 in Duffie et al. \cite{DufFilSch}.
This way of proving is somewhat complicated that's why we decided to find another way presented above.

For all \ $n,p\in\ZZ_+$, \ using the independence of \ $L$ \ and \ $B$, \ by It\^{o}'s formula, we have
 \begin{align*}
  \dd(Y_t^n X_t^p)
    & = nY_t^{n-1} X_t^p \big( (a-bY_t)\,\dd t + \sqrt{Y_t}\,\dd L_t\big)
       + p Y_t^n X_t^{p-1} \big( (m-\theta X_t)\,\dd t + \!\sqrt{Y_t}\,\dd B_t \big)\\
    &\phantom{=\;}
       + \frac{n(n-1)}{2} Y_t^{n-2} X_t^p Y_t\,\dd t
       + \frac{p(p-1)}{2} Y_t^n X_t^{p-2} Y_t\,\dd t \\
    & =\Big(nY_t^{n-1} (a-bY_t) X_t^p
            + pY_t^n(m-\theta X_t)X_t^{p-1}
            + \frac{n(n-1)}{2}Y_t^{n-1} X_t^p\\
    &\phantom{=\Big(\,}
        + \frac{p(p-1)}{2}Y_t^{n+1} X_t^{p-2}
        \Big)\dd t + nY_t^{n-1/2}X_t^p\,\dd L_t
       + pY_t^{n+1/2}X_t^{p-1}\,\dd B_t
 \end{align*}
 for \ $t\geq 0$.
\ Writing the SDE above in an integrated form and taking expectations of both of its sides, we have
 \begin{align*}
   \EE(Y_t^n X_t^p)
     - \EE(Y_0^n X_0^p)
    &  =  \int_0^t \Big[
           an \EE(Y_u^{n-1}X_u^p)
           -bn \EE(Y_u^n X_u^p)
           + pm \EE(Y_u^nX_u^{p-1})
           - p\theta \EE(Y_u^nX_u^p)\\
    &\phantom{=  \int_0^t \Big[\;}
          + \frac{n(n-1)}{2}\EE(Y_u^{n-1} X_u^p)
          + \frac{p(p-1)}{2}\EE(Y_u^{n+1} X_u^{p-2}) \Big]\,\dd u,
          \qquad t\geq 0,
 \end{align*}
 where we used that
 \[
   \left(\int_0^t Y_u^{n-1/2}X_u^p\,\dd L_u\right)_{t\geq 0}
    \qquad \text{and} \qquad
   \left(\int_0^t Y_u^{n+1/2}X_u^{p-1}\,\dd B_u\right)_{t\geq 0}
 \]
 are continuous square integrable martingales due to \eqref{help65}, see, e.g., Ikeda and Watanabe \cite[page 55]{IkeWat}.
Introduce the functions \ $f_{n,p}(t) := \EE(Y_t^n X_t^p)$, \ $t\in\RR_+$, \ for
 \ $n,p\in\ZZ_+$.
\ Then we have
 \begin{align*}
   f_{n,p}'(t)
        = - (bn+p\theta) f_{n,p}(t)
          + \left( an + \frac{n(n-1)}{2} \right) f_{n-1,p}(t)
          + pm f_{n,p-1}(t)
          + \frac{p(p-1)}{2} f_{n+1,p-2}(t)
 \end{align*}
 for \ $t \in \RR_+$, \ where \ $f_{k,\ell}(t) := 0$ \ if \ $k,\ell\in\ZZ$ \ with
 \ $k < 0$ \ or \ $\ell < 0$.
\ Hence for all \ $M \in \NN$, \ the functions \ $f_{n,p}$, \ $n,p\in\ZZ_+$
 \ with \ $n+p\leq M$ \ satisfy a homogeneous linear system of differential
 equations with constant coefficients.
For example, if \ $M=2$ \ then
 \[
   \begin{bmatrix}
    f_{0,0}'(t) \\ f_{1,0}'(t) \\ f_{0,1}'(t) \\
    f_{2,0}'(t) \\ f_{1,1}'(t) \\ f_{0,2}'(t)
   \end{bmatrix}
   =\begin{bmatrix}
     0 & 0 & 0 & 0 & 0 & 0 \\
     a & -b & 0 & 0 & 0 & 0 \\
     m & 0 & -\theta & 0 & 0 & 0 \\
     0 & 2a+1 & 0 & -2b & 0 & 0 \\
     0 & m & a & 0 & -b-\theta & 0 \\
     0 & 1 & 2m & 0 & 0 & -2\theta
    \end{bmatrix}
    \begin{bmatrix}
     f_{0,0}(t) \\ f_{1,0}(t) \\ f_{0,1}(t) \\
     f_{2,0}(t) \\ f_{1,1}(t) \\ f_{0,2}(t)
    \end{bmatrix} , \qquad t \in \RR_+ .
 \]
Thus, for all \ $n,p\in \ZZ_+$, \ the function \ $f_{n,p}$ \ is a linear
 combination of the functions \ $\ee^{-(kb+\ell\theta)t}$, \ $t\in \RR_+$,
 \ $k,\ell\in\ZZ_+$ \ with \ $k+\ell\leq n+p$, \ since the eigenvalues of the
 coefficient matrix of the above mentioned system of differential equations are
 \ $-(kb+\ell\theta)$, \ $k,\ell\in\ZZ_+$ \ with \ $k+\ell\leq M$.
\ Consequently, for all \ $n,p\in \ZZ_+$, \ the function \ $f_{n,p}$ \ is
 bounded and the limit \ $\lim_{t\to\infty} f_{n,p}(t)$ \ exists
 and finite.
By the moment convergence theorem
 (see, e.g., Stroock \cite[Lemma 2.2.1]{STR}),
 \ $\lim_{t\to\infty} f_{n,p}(t)
    =\lim_{t\to\infty} \EE(Y_t^n X_t^p)
    = \EE(Y_\infty^n X_\infty^p)$, \ $n,p\in\ZZ_+$.
\ Indeed, by Theorem \ref{Thm_ergodic1} and the continuous mapping theorem,
 \ $Y_t^n X_t^p \distr Y_\infty^n X_\infty^p$ \ as \ $t\to\infty$, \ and
 the family \ $\{Y_t^n X_t^p : t\in\RR_+\}$ \ is uniformly integrable.
This latter fact follows from the boundedness of the function \ $f_{2n,2p}$,
 \ see, e.g., Stroock \cite[condition (2.2.5)]{STR}.
Hence we conclude that all the mixed moments of \ $(Y_\infty, X_\infty)$ \ are finite.

Finally, we calculate the moments listed in the theorem.
Let us consider again the 2-dimensional affine diffusion model \eqref{2dim_affine} with \ $\alpha=2$,
 \ $a>0$, \ $b>0$, $m\in\RR$, \ $\theta>0$, \ and with a random initial value
 \ $(Y_0,X_0)$ \ independent of \ $(L_t,B_t)_{t\geq 0}$ \ having the same distribution as that of \ $(Y_\infty,X_\infty)$.
\ Then, by Theorem \ref{Thm_ergodic1}, the process  \ $(Y_t,X_t)_{t\geq 0}$ \ is
 strictly stationary, and hence, \ $f_{n,p}(t) = \EE(Y_\infty^n X_\infty^p)$,
 \ for all \ $t\in\RR_+$ \ and \ $n,p\in\ZZ_+$.
\ The above system of differential equations for the functions
 \ $f_{n,p}$, \ $n,p\in\ZZ_+$, \ yields
 \begin{align}\label{help51}
  \begin{split}
   \EE(Y_\infty^n X_\infty^p)
   & = \frac{1}{bn+p\theta}
       \Bigg( \left( an + \frac{n(n-1)}{2}\right) \EE(Y_\infty^{n-1}X_\infty^p)
              + mp \EE(Y_\infty^nX_\infty^{p-1}) \\
   &\phantom{= \frac{1}{bn+p\theta} \Bigg(\;}
              + \frac{p(p-1)}{2}\EE(Y_\infty^{n+1} X_\infty^{p-2})
       \Bigg)
 \end{split}
 \end{align}
 for all \ $n,p\in\ZZ_+$.
\ By \eqref{help51}, one can calculate the moments listed in the theorem.

Finally, we note that for calculating the moments \ $\EE(Y_\infty^n X_\infty^p)<\infty$, \ $n,p\in\ZZ_+$, \
 we could have used formula (4.4) in Filipovi\'{c} et al. \cite{FilMaySch} which gives a formal representation
 of the polynomial moments of \ $(Y_t,X_t)$, \ $t\in\RR_+$.
\ The idea behind this formal representation is that the infinitesimal generator of the
 affine process \ $(Y,X)$ \ formally maps the finite-dimensional linear space of all polynomials in \ $(y,x)\in\RR_+\times\RR$
 \ of degree less than or equal to \ $k$ \ into itself, where \ $k\in\NN$.
\ For a more general class of time-homogeneous Markov processes having this property,
 for the so-called polynomial processes, see Cuchiero et al. \cite{CucKelTei}.

We also remark that the moments of \ $Y_\infty$ \ could have been calculated directly using that
 \ $Y_\infty$ \ has Gamma distribution with parameters \ $2a$ \ and \ $2b$.
\proofend

\section*{Acknowledgements}
We are undoubtedly grateful for the referee for pointing out a mistake in the proof of Theorem \ref{Thm_ergodic2},
 and also for his/her several valuable comments that have led to an improvement of the manuscript.

\end{document}